\documentclass[journal,onecolumn]{IEEEtran}
\usepackage{amssymb}
\usepackage{graphicx}
\usepackage{epsfig}
\usepackage{subcaption}
\usepackage{multicol, blindtext}
\usepackage{lipsum}
\usepackage{multirow}
\usepackage{amsmath}
\usepackage{amsfonts}
\usepackage{amssymb}
\usepackage{hyperref}
\usepackage{bm}
\usepackage[usenames]{color}
\newtheorem{example}{Example}
\newtheorem{definition}{Definition}
\newtheorem{proposition}{Proposition}
\newtheorem{lemma}{Lemma}
\newtheorem{theorem}{Theorem}
\newtheorem{corollary}{Corollary}
\usepackage{algorithm,algorithmic}
\usepackage{float}
\floatstyle{boxed}
\restylefloat{figure}
\renewcommand{\algorithmicrequire}{\textbf{Input: }}
\renewcommand{\algorithmicensure}{\textbf{Output: }}

   \newtheorem{assumption}{Assumption}

\newcommand{\oomit}[1]{}

\begin{document}

\title{Inner-Approximating Reachable Sets for Polynomial Systems with Time-Varying Uncertainties}

\author{
Bai Xue$^{1,2}$ and Martin Fr\"anzle$^{3}$ and Naijun Zhan$^{1,2}$
\thanks{1. State Key Lab. of Computer Science, Institute of Software, CAS, China.
        Email: \{xuebai,znj\}@ios.ac.cn }
\thanks{2. University of Chinese Academy of Sciences, Beijing, China.}
\thanks{3. Carl von Ossietzky Universit\"at Oldenburg, Germany.
     Email: fraenzle@informatik.uni-oldenburg.de}
}

\maketitle

\IEEEpeerreviewmaketitle

\newcommand\mf[1]{\textcolor{blue}{\footnotesize MF: #1}}

\begin{abstract}
In this paper we propose a convex programming based method to address a long-standing problem of inner-approximating backward reachable sets of state-constrained polynomial systems subject to time-varying uncertainties. The backward reachable set is a set of states, from which all trajectories starting will surely enter a target region at the end of a given time horizon without violating a set of state constraints in spite of the actions of uncertainties. It is equal to  the zero sub-level set of the unique Lipschitz viscosity solution to a Hamilton-Jacobi partial differential equation (HJE). We show that inner-approximations of the backward reachable set can be formed by zero sub-level sets of its viscosity super-solutions. Consequently, we reduce the inner-approximation problem to a problem of synthesizing polynomial viscosity super-solutions to this HJE. Such a polynomial solution in our method is synthesized by solving a single semi-definite program. We also prove that polynomial solutions to the formulated semi-definite program exist and can produce a convergent sequence of inner-approximations to the interior of the backward reachable set in measure under appropriate assumptions. This is the main contribution of this work. Several illustrative examples demonstrate the merits of our approach.
\end{abstract}

\begin{IEEEkeywords}
Reachability Analysis; Uncertainties; Convex Computations
\end{IEEEkeywords}

\section{Introduction}
\label{introduction}
Reachability analysis, which derives verdicts about the states reachable in a dynamical system, has received growing interest in recent years. It has many applications in engineering problems, especially concerning safety-critical systems including aeronautics, automotive, medical devices and industrial process control 
\cite{lunze2009}. Consequently, attention from scientists across multiple disciplines has been devoted to the problem of performing reachability analysis. Performing outer- and inner-approximate reachability analysis is an enabler for detecting whether the system of interest will always avoid unsafe states when started from a specified set of initial states or whether it satisfies a temporal-logic formula \cite{eggers2015}, as well as for computing the set of initial configurations that reach desired configurations while respecting a set of constraints \cite{aubin2011}. The former is generally referred to as  the safety verification problem, which has traditionally attracted more attention. As a result, significant advances of outer-approximate reachability analysis techniques for both linear and nonlinear systems have been reported in the literature based on various representations of sets such as intervals \cite{Ramdani2009}, zonotopes \cite{Althoff08}, polyhedra and support functions for polyhedral sets \cite{Dang10,Frehse15}, ellipsoids \cite{Kurzhanski00}, level sets \cite{mitchell2005}, Taylor models \cite{Xin12} and semi-algebraic sets \cite{Wang2013,henrion2014}. Computational methods for inner-approximations have received increasing attention just recently, e.g., \cite{Wang2013,Korda13,Goubault14,Xin2014,Xue16}. It nevertheless has a wide range of practical applications including collision avoidance and surveillance. However, the development of numerical tools, which tractably inner-approximate the reachable set for state-constrained systems with time-varying uncertainties, has been challenging and is still an open area of research.

Besides, in real physical world physical systems often have certain level of desired performances. Unfortunately, it is demanding to model their dynamics exactly due to physical limitations such as imperfections in sensing equipment and incomplete information, especially in
fluctuating environments. Consequently, engineering designs based on abstracted mathematical models without taking these uncertainties into account may lead to incorrect operations of physical systems. Abstracting these uncertainties as time-varying parameters (e.g., \cite{rocca2018}) and incorporating them into the model is a popular means to compensate for the inability to construct exact models.

In this paper we focus our attention on inner-approximating backward reachable sets for state-constrained polynomial systems  with time-varying uncertainties. The backward reachable set is the set of states such that trajectories originating from it surely hit a  target region after a specified time duration without violating a  set of state constraints in spite of the actions of the uncertainties. Such sets are particularly useful to identify decisions that are ``robust'' against noise parameters. In order to compute the backward reachable set, in this paper we first make use of Kirszbraun's extension theorem for Lipschitz maps to characterize the backward reachable set as the zero sub-level set of the unique Lipschitz viscosity solution to a HJE. Such HJE could be regarded as a special case of the HJE in \cite{margellos2011} considering competing inputs (uncertainty and control) and time-invariant state constraints. Since it is nontrivial, even impossible to find the viscosity solution, we then propose a novel semi-definite programming based method to compute its polynomial viscosity super-solutions, whose zero sub-level sets form inner approximations of the backward reachable set. An inner-approximation of the backward reachable set in our method can be obtained by solving a single semi-definite program consisting of linear matrix inequalities. Compared to traditionally grid-based numerical methods, the benefits of our method are overall the convexity of the problem of finding the backward reachable set. We further prove that polynomial solutions to the formulated semi-definite program exist and can generate a convergent sequence of inner-approximations to the interior of the backward reachable set in measure under appropriate assumptions. This is the main contribution of this work. Finally, several illustrative examples evaluate the performance and the merits of our approach. 

\subsection*{Related Work}
As mentioned above, inner-approximate reachability analysis of ordinary differential equations subject to time-varying uncertainties and state constraints, is still in its infancy and thus provides an open area of research.

For ordinary differential equations free of time-varying uncertainties and state constraints, \cite{Goubault14,goubault2017} proposed a method based on modal intervals with affine forms to inner-approximate reachable sets using intervals. By making use of the homeomorphism property of the solution mapping, a boundary based reachability analysis method was proposed to inner-approximate reachable sets with polytopes in \cite{Xue16}, and it was extended to a class of delay differential equations in  \cite{xue2017safe}. Since  reachable sets of nonlinear systems tend to be non-convex, the above mentioned methods based on convex set representations may result in poor approximations. As accuracy is also an important factor in performing reachability analysis (e.g.,\cite{schupp2015,maler2014}), more complex shapes of representations such as  Taylor models and semi-algebraic sets are desirable. \cite{Xin2014} proposed a Taylor model backward flowpipe method to compute inner-approximations. \cite{Wang2013} proposed an iterative method, with each iteration relying on solving semi-definite programming problems, to compute semi-algebraic inner-approximations for polynomial systems using the advection map of the given dynamical system. \cite{xue2017} extended the method in \cite{Xue16} to compute semi-algebraic inner-approximations of reachable sets for polynomial systems and beyond by solving semi-definite programming problems. Recently, \cite{BMN18} formulated the problem of solving HJEs as a semi-definite program to compute inner-approximations for polynomial systems. For state-constrained polynomial systems without time-varying uncertainties, \cite{Korda13} computed inner-approximations of the region of attraction to a target set by solving semi-definite programs. In contrast to the aforementioned approaches, our approach in this paper targets state-constrained systems subject to time-varying uncertainties.

The reachability analysis for state-constrained nonlinear systems with time-varying uncertainties is more challenging. An attractive way to address this problem is by formulating reachable sets as sub-level sets of viscosity solutions of HJEs, e.g, \cite{mitchell2005,bokanowski2010,margellos2011,aubin2011,fisac2015,WZF19}. The Hamilton-Jacobi reachability methods are capable of dealing with general nonlinear systems with state constraints and competing inputs. However, existing numerical methods for addressing HJEs generally require gridding the state space and hence their time and memory complexity grow exponentially with the state dimension. Our approach in this paper tackles the finite time horizon reachability problem of state-constrained polynomial systems with time-varying uncertainties. Rather than solving HJEs directly, our approach reformulates the problem of solving HJEs as a semi-definite programming problem, which falls within the convex programming framework and can be efficiently solved by interior-point methods in polynomial time. Polynomial solutions to the formulated semi-definite programs exist and can produce a convergent sequence of inner-approximations to the interior of the backward reachable set in measure under appropriate assumptions. Recently, based on a derived HJE, \cite{WZF19} proposed a semi-definite programming based method to compute inner-approximations of the maximal robust invariant set over the infinite time horizon for state-constrained polynomial systems with time-varying uncertainties. However, the existence of polynomial solutions to the constructed semi-definite program in \cite{WZF19} is not guaranteed.

Another area that is relevance to the topic of this paper is the computation of regions of attraction for systems subject to uncertainties \cite{chesi2004,topcu2007,topcu2010,chesi2013}. These methods rely upon the generation or evaluation of pre-constructed Lyapunov functions to compute inner-approximations of the region of attraction over the infinite time horizon. This requires checking Lyapunov's criteria for polynomial systems by using sum-of-squares programming, which results in a bilinear optimization problem that is usually solved using some form of alteration, e.g., \cite{topcu2010}. These sum-of-squares programming based methods suffer from the same issue as that in \cite{WZF19}. The existence of polynomial solutions to the constructed sum-of-squares programming is not guaranteed. 

This paper is structured as follows. The reachability problem of interest is formally stated in Section \ref{pre}, and then formulated within the Hamilton-Jacobi reachability framework in Section \ref{theroy and algorithm}. In Section \ref{IACR} we show that the interior of the backward reachable set can be approximated from inside in measure by a sequence of zero sub-level sets of solutions to a semi-definite program under appropriate assumptions. After demonstrating our approach on several illustrative examples in Section \ref{examples}, we conclude our paper in Section \ref{conclusion}.

\section{Preliminaries}
\label{pre}
\subsection{System Dynamics}
 In this section we mainly present an introduction to backward reachable sets. The following notation will be used throughout this paper: For a set $\Delta$, $\partial \Delta$ denotes its boundary. $R_k[\cdot]$ represents the set of real 
polynomials of total degree $\leq k$ in variables given by the argument. The symbol $\mathbb{R}[\cdot]$ denotes the ring of polynomials in variables given by the argument. $\mathbb{N}$ denotes the set of nonnegative integers. The space of continuously differentiable functions on a set $X$ is denoted by $C^{\infty}(X)$. The difference of two sets of $A$ and $B$ is denoted by $A\setminus B$. $\mu(A)$ denotes the Lebesgue measure on $A\subset \mathbb{R}^n$. Vectors are denoted by boldface letters.

In this paper we consider the following system:
\begin{equation}
\label{systems}
\dot{\bm{x}}(s)=\bm{f}(\bm{x}(s),\bm{d}(s)),a.e., s\in [0,T],
\end{equation}
where for each $s\in [0,T]$, $\bm{x}(s)\in \mathcal{X}_s$ and $\bm{d}(s)\in \mathcal{D}$, $\mathcal{X}_s$ and $\mathcal{D}$ are respectively compact subsets of $\mathbb{R}^n$ and $\mathbb{R}^m$ for some positive integers $n$ and $m$.

We assume that each entry of the vector field $\bm{f}$ is polynomial, i.e., $f_i\in \mathbb{R}[\bm{x},\bm{d}], i=1,\ldots,n$. It is evident that the map $\bm{f}$ satisfies the following two properties:
\begin{enumerate}
\item $\bm{f}$ is continuous;
\item $\bm{f}$ is locally Lipschitz on $\bm{x}$ uniformly on $\bm{d}$, that is, for each compact subset $\mathcal{X}$ of $\mathbb{R}^n$ there is some constant $L$ such that \[\|\bm{f}(\bm{x},\bm{d})-\bm{f}(\bm{z},\bm{d})\|\leq L\|\bm{x}-\bm{z}\|,\forall \bm{x},\bm{z}\in \mathcal{X}, \forall\bm{d}\in \mathcal{D},\] where $\|\cdot\|$ denotes the usual Euclidean norm.
\end{enumerate}

For $t\in [0,T]$, the time-varying state and uncertainty constraint sets $\mathcal{X}_t$ and $\mathcal{D}$ are basic compact semi-algebraic sets, i.e. 
\[\mathcal{X}_t:=\{\bm{x}\in \mathbb{R}^n\mid g_i(\bm{x},t)\leq 0, i=1,\ldots, n_{\mathcal{X}}\}\]
\[\mathcal{D}:=\{\bm{d}\in \mathbb{R}^m\mid h_i(\bm{d})\geq 0,i=1,\ldots,n_{\mathcal{D}}\}\]
with $g_i\in \mathbb{R}[\bm{x},t]$ and $h_i\in \mathbb{R}[\bm{d}]$. Also, $\partial \mathcal{X}_t=\cup_{i=1}^{n_{\mathcal{X}}} \{\bm{x}\in \mathcal{X}_t\mid g_i(\bm{x},t)=0\}$. The terminal state $\bm{x}(T)$ is constrained to lie in the basic semi-algebraic set $\mathtt{TR}$, where \[\mathtt{TR}:=\{\bm{x}\in \mathbb{R}^n\mid l_i(\bm{x})\leq 0, i=1,\ldots,n_{\mathtt{TR}}\}\]
with $l_i\in \mathbb{R}[\bm{x}]$ and $\partial \mathtt{TR}=\cup_{i=1}^{n_{\mathtt{TR}}} \{\bm{x}\in \mathtt{TR}\mid l_i(\bm{x})=0\}$.

Let $\mathcal{M}_{t}$ be the set of measurable functions $\bm{d}:[t,T] \mapsto \mathcal{D}$, where $T>t$. We will call functions $\bm{d}\in \mathcal{M}_{t}$ time-varying uncertainties. For each $\bm{d}\in \mathcal{M}_{t}$, we denote by $\bm{y}_{\bm{x},t}^{\bm{d}}(s)$ the solution at time $s\in [t,T]$ of \eqref{systems} starting from the state $\bm{x}$ at time $t$.

The problem we attempt to address is to compute the backward reachable set $\mathcal{R}_0$ such that all trajectories starting from it at time $t=0$ will enter the target region $\mathtt{TR}$ after the time duration of $T$ while staying inside the set $\mathcal{X}_s$ for $s\in [0,T]$, despite the actions of uncertainties.
\begin{definition}
\label{BRS}
The backward reachable set $\mathcal{R}_0$ of the target region $\mathtt{TR}$ at time $t=0$ is presented as follows:
\begin{equation}
\begin{split}
\mathcal{R}_0:=\{\bm{x}_0|\forall \bm{d}\in \mathcal{M}_{0},&\bm{y}_{\bm{x}_0,0}^{\bm{d}}(0)=\bm{x}_0, \bm{y}_{\bm{x}_0,0}^{\bm{d}}(T)\in \mathtt{TR},\bm{y}_{\bm{x}_0,0}^{\bm{d}}(s)\in \mathcal{X}_s~\text{for}~s\in [0,T]\}.
\end{split}
\end{equation}
\end{definition}

The backward reachable set in Definition \ref{BRS} differs from the constrained controlled region of attraction $\mathcal{R}'_0$ in \cite{henrion2014,pauwels2017}. The constrained controlled region of attraction $\mathcal{R}'_0$ is the set of initial states that can be driven with an admissible control to a specified target set without leaving the state-constrained set, i.e. 
\begin{equation}
\begin{split}
\mathcal{R}'_0=\{\bm{x}_0|\exists \bm{d}\in \mathcal{M}_{0},&\bm{\phi}_{\bm{x}_0,0}^{\bm{d}}(0)=\bm{x}_0, \bm{\phi}_{\bm{x}_0,0}^{\bm{d}}(T)\in \mathtt{TR},\bm{\phi}_{\bm{x}_0,0}^{\bm{d}}(s)\in \mathcal{X}_s~\text{for}~s\in [0,T]\}.
\end{split}
\end{equation}
Obviously, $\mathcal{R}_0\subset \mathcal{R}'_0$. In \cite{henrion2014,pauwels2017}, an outer approximation of the set  $\mathcal{R}'_0$ is computed by solving a single semi-definite program, which is constructed from occupation measure.

It is in general impossible to obtain the backward reachable set $\mathcal{R}_0$ since an appropriate closed-form solution to \eqref{systems} may not be available. We therefore resort to the computation of an inner approximation of the backward reachable set. We opt for inner approximations as they preserve the desired property of the backward reachable set, namely that all possible trajectories starting from them enter $\mathtt{TR}$ after the time duration of $T$ while not leaving the set $\mathcal{X}_t$ for $t\in [0,T]$.

\section{Hamilton-Jacobi Type Equations}
\label{theroy and algorithm}
In this section we mainly introduce the reformulation of the backward reachable set $\mathcal{R}_0$ as the zero sub-level set of the viscosity solution to a Hamilton-Jacobi type partial differential equation. Like in \cite{WZF19}, we use Kirszbraun's extension theorem to characterize the backward reachable set $\mathcal{R}_0$ as the zero sub-level set of the unique Lipschitz viscosity solution to a HJE. 

As $\bm{f}\in \mathbb{R}[\bm{x},\bm{d}]$ in system \eqref{systems}, $\bm{f}$ is locally Lipschitz continuous over the state variable $\bm{x}$. Therefore, the global solution $\bm{\phi}_{\bm{x}_0}^{d}(t)$ over $t\in [0,\infty)$ to system \eqref{systems} is not guaranteed to exist for any initial state $\bm{x}_0\in \mathbb{R}^n$. This hinders the construction of Hamilton-Jacobi equations. To address this issue, we first construct an auxiliary vector field $\bm{F}(\bm{x},\bm{d}):\mathbb{R}^n\times \mathcal{D}\mapsto \mathbb{R}^n$, which is globally Lipschitz on $\bm{x}\in \mathbb{R}^n$ uniformly on $\bm{d}\in \mathcal{D}$, i.e. there exists a constant $L_F$ such that
 \[\|\bm{F}(\bm{x}_1,\bm{d})-\bm{F}(\bm{x}_2,\bm{d})\|\leq L_F\|\bm{x}_1-\bm{x}_2\|, \forall \bm{x}_1,\bm{x}_2\in \mathbb{R}^n, \forall \bm{d}\in \mathcal{D},\]
where $\|\cdot\|$ denotes the usual Euclidean norm. Moreover, the trajectories governed by $\dot{\bm{x}}(s)=\bm{F}(\bm{x}(s),\bm{d}(s))$ coincide with the trajectories generated by $\dot{\bm{x}}(s)=\bm{f}(\bm{x}(s),\bm{d}(s))$ over a local state space. Thanks to Kirszbraun's theorem \cite{fremlin2011}, which is stated in Theorem \ref{kri}, the existence of such function is ensured.
 \begin{theorem}[Kirszbraun's Theorem]
 \label{kri}
 Let $H_1$ and $H_2$ be Hilbert spaces, $A\subset H_1$ a set and $\bm{f}':A\mapsto H_2$ a function. Suppose that $\gamma\geq 0$ is such that $\|\bm{f}'(\bm{x})-\bm{f}'(\bm{y})\|\leq \gamma \|\bm{x}-\bm{y}\|$ for $\bm{x},\bm{y}\in A$. Then there is a function $\bm{F}':H_1\mapsto H_2$ such that $\bm{F}'(\bm{x})=\bm{f}'(\bm{x})$ for $\bm{x}\in A$ and $\|\bm{F}'(\bm{x})-\bm{F}'(\bm{y})\|\leq \gamma \|\bm{x}-\bm{y}\|$ for all $\bm{x},\bm{y}\in H_1$.
 \end{theorem}

 Thus, rather than considering system \eqref{systems}, in this subsection we take into account an auxiliary system:
 \begin{equation}
\label{sys1}
\dot{\bm{x}}(s)=\bm{F}(\bm{x}(s),\bm{d}(s)), a.e., s\in [0,T],
 \end{equation}
 where for each $s\in [0,T]$, $\bm{x}(s)\in \mathcal{X}_s$ and $\bm{d}(s) \in \mathcal{D}$, and where $\mathcal{X}_s$ and $\mathcal{D}$ are respectively compact subsets of $\mathbb{R}^n$ and $\mathbb{R}^m$. The map $\bm{F}:\mathbb{R}^n\times \mathcal{D} \mapsto \mathbb{R}^n$ is assumed to satisfy the following three properties:
 \begin{enumerate}
 \item $\bm{F}: \mathbb{R}^n\times \mathcal{D}\mapsto \mathbb{R}^n$ is continuous;
 \item $\bm{F}$ is globally Lipschitz continuous on $\bm{x}\in \mathbb{R}^n$ uniformly on $\bm{d} \in \mathcal{D}$, that is, there is some constant
 $L_{F}$ such that
 \[\|\bm{F}(\bm{x},\bm{d})-\bm{F}(\bm{y},\bm{d})\|\leq L_F\|\bm{x}-\bm{y}\|\]
for all $\bm{x},\bm{y}\in \mathbb{R}^n$ and all $\bm{d}\in \mathcal{D}$, where $\|\cdot\|$ denotes the usual Euclidean norm;
\item $\bm{F}(\bm{x},\bm{d})=\bm{f}(\bm{x},\bm{d})$ over $\bm{x}\in B(\bm{0},R)$ and $\bm{d}\in \mathcal{D}$, where
\begin{equation}
\label{B}
B(\bm{0},R)=\{\bm{x}\in \mathbb{R}^{n}\mid g_R(\bm{x})\geq 0\},
\end{equation}
thereof, $g_R(\bm{x})=R-\sum_{i=1}^nx_i^2$ and $R$ is a positive number such that $\mathcal{X}_t\subseteq B(\bm{0},R)$ for $t\in [0,T]$ with $\partial \mathcal{X}_t\cap \partial B(\bm{0},R)=\emptyset$. 
Note that $R$ exists since $\cup_{t\in [0,T]}\mathcal{X}_t$ is compact due to Lemma 1 in \cite{fisac2015}.
 \end{enumerate}

The set $B(\bm{0},R)$ in \eqref{B} plays three important roles in our approach.
\begin{enumerate}
\item The condition $\mathcal{X}_t\subseteq B(\bm{0},R)$ for $t\in [0,T]$ guarantees that the backward reachable set $\mathcal{R}_0$ for system \eqref{systems} can be characterized by trajectories to the auxiliary system \eqref{sys1}, as formulated in Proposition \ref{relation}.
\item The condition $\partial \mathcal{X}_t\cap \partial B(\bm{0},R)=\emptyset$ for $t\in [0,T]$ assures that the zero sub-level set of a polynomial, which is computed by solving \eqref{sos} in Subsection \ref{SSP}, is an inner-approximation of the backward reachable set $\mathcal{R}_0$, as stated in Corollary \ref{inner2} in Subsection \ref{SSP}.
\item The condition $B(\bm{0},R)=\{\bm{x}\in \mathbb{R}^{n}\mid g_R(\bm{x})\geq 0\}$ with $g_R(\bm{x})=R-\sum_{i=1}^nx_i^2$ is used to guarantee the existence of solutions to the semi-definite program \eqref{sos} in Subsection \ref{SSP}. It is useful in justifying Corollary \ref{unifrom} in Subsection \ref{CA}.
\end{enumerate}

Now we know that for any $\bm{d}\in \mathcal{D}$ and any $\bm{x}\in \mathbb{R}^n$, there exists a unique absolutely continuous trajectory $\bm{y}(s)=\bm{\phi}_{\bm{x},t}^{\bm{d}}(s)$ satisfying \eqref{sys1} for almost all $s\geq t$ and $\bm{y}(t)=\bm{x}.$
\begin{definition}
\label{re}
For $T>t$ with $t\geq 0$, the set of states, which are visited by trajectories on $[t,T]$ starting from $\bm{x}$, is denoted as
$S_{[t,T]}(\bm{x}):=\{\bm{y}\in \mathbb{R}^n\mid \bm{y}=\bm{\phi}_{\bm{x},t}^{\bm{d}}(s)$ is absolutely continuous, satisfies \eqref{sys1} for some $\bm{d}\in \mathcal{M}_{t}$, $\bm{y}(\bm{t})=\bm{x}$ \text{~and~} $s\in [t,T]\}$.
\end{definition}
 
 Again under above assumption, for $T>0$ and $\bm{x}\in \mathbb{R}^n$, $S_{[0,T]}(\bm{x})$ is a compact set in the Soblov space $W^{1,1}(0,T)$ for the topology of $C([0,T];\mathbb{R}^n)$.

Next, consider the backward reachable set $\mathcal{R}_t$ of $\mathtt{TR}$ at time $t$ for system \eqref{systems}, which is the set of states such that all trajectories starting from  it at time $t$ will enter the target region $\mathtt{TR}$ after the time duration of $T-t$ while not leaving the state constraint set $\mathcal{X}_{s}$ for $ s\in [t,T]$, i.e.,
\begin{equation}
\label{max}
\begin{split}
\mathcal{R}_t:=\{\bm{x}\in \mathbb{R}^n|& \forall \bm{d}\in \mathcal{M}_t, \forall s\in [t,T], \bm{y}_{\bm{x},t}^{\bm{d}}(s)\in \mathcal{X}_s, \bm{y}_{\bm{x},t}^{\bm{d}}(T)\in \mathtt{TR}\},
\end{split}
\end{equation}
where $\bm{y}_{\bm{x},t}^{\bm{d}}(\cdot):[t,T]\mapsto \mathbb{R}^n$ is the solution to system \eqref{systems} with  $\bm{d}\in \mathcal{M}_t$ for the time interval $[t,T]$.

Let's present a value function $u(\bm{x},t)$ defined below,
\begin{equation}
\label{u}
\begin{split}
u(\bm{x},t):=\sup_{\bm{d}\in \mathcal{M}_t}&\max\big\{\max_{i\in \{1,\ldots,n_{\mathtt{TR}}\}}\{l_i(\bm{\phi}_{\bm{x},t}^{\bm{d}}(T))\},\max_{i\in \{1,\ldots,n_{\mathcal{X}}\}} \{\max_{s \in [t,T]}g_i(\bm{\phi}_{\bm{x},t}^{\bm{d}}(s),s)\}\big\}.
\end{split}
\end{equation}

Proposition \ref{relation} builds a relationship between the value function $u(\bm{x},t)$ and the backward reachable set $\mathcal{R}_t$.
\begin{proposition}
\label{relation}
$\mathcal{R}_t=\{\bm{x}\in \mathbb{R}^n\mid u(\bm{x},t)\leq 0\}$.
\end{proposition}
\begin{IEEEproof}
Obviously, according to the relationship between \eqref{systems} and \eqref{sys1}, if the trajectory $\bm{\phi}_{\bm{x},t}^{\bm{d}}(s)$ to system \eqref{systems} stays in the set $B(\bm{0},R)$ for $ s \in [t,T]$, where $\bm{d}\in \mathcal{M}_t$, we have $\bm{\phi}_{\bm{x},t}^{\bm{d}}(s)=\bm{y}_{\bm{x},t}^{\bm{d}}(s)$. Also, since $\mathcal{X}_\tau \subset B(\bm{0},R)$ for $ \tau\in [t,T]$ and $\mathtt{TR}\subset B(\bm{0},R)$, we have $$\{\bm{x}\in \mathbb{R}^n \mid u(\bm{x},t)\leq 0\}=\{\bm{x}\in \mathbb{R}^n\mid v(\bm{x},t)\leq 0\},$$ where
\begin{equation}
\label{1v}
\begin{split}
v(\bm{x},t)=\sup_{\bm{d}\in \mathcal{M}_t}&\max\big\{\max_{i\in \{1,\ldots,n_{\mathtt{TR}}\}}\{l_i(\bm{y}_{\bm{x},t}^{\bm{d}}(T))\},\max_{i\in \{1,\ldots,n_{\mathcal{X}}\}} \{\max_{s \in [t,T]}g_i(\bm{y}_{\bm{x},t}^{\bm{d}}(s),s)\}\big\}.
\end{split}
\end{equation}
Thus, it is enough to prove that $\mathcal{R}_t=\{\bm{x}\mid v(\bm{x},t)\leq 0\}$.

If $\bm{x}\in \mathcal{R}_t$, according to the definition of $\mathcal{R}_t$, i.e. \eqref{max}, we can deduce that for all $\bm{d}\in \mathcal{M}_t$,
\begin{equation}
\begin{split}
\max\big\{&\max_{i\in \{1,\ldots,n_{\mathtt{TR}}\}}\{l_i(\bm{y}_{\bm{x},t}^{\bm{d}}(T))\},\max_{i\in \{1,\ldots,n_{\mathcal{X}}\}}\{\max_{s\in [t,T]}g_i(\bm{y}_{\bm{x},t}^{\bm{d}}(s),s)\}\big\}\leq 0,
\end{split}
\end{equation}
implying that
\begin{equation}
\label{max1}
\begin{split}
\sup_{\bm{d}\in \mathcal{M}_t}\max\big\{&\max_{i\in \{1,\ldots,n_{\mathtt{TR}}\}}\{l_i(\bm{y}_{\bm{x},t}^{\bm{d}}(T))\},\max_{i\in \{1,\ldots,n_{\mathcal{X}}\}}\{\max_{s\in [t,T]}g_{i}(\bm{y}_{\bm{x},t}^{\bm{d}}(s),s)\}\big\}\leq 0.
\end{split}
\end{equation}
Consequently, $v(\bm{x},t)\leq 0$.

On the other hand, if $\bm{x}\in \{\bm{x}\in \mathbb{R}^n| v(\bm{x},t)\leq 0\}$,  according to \eqref{1v}, then $$\max_{i\in \{1,\ldots,n_{\mathtt{TR}}\}}\{l_i(\bm{y}_{\bm{x},t}^{\bm{d}}(T))\}\leq 0, \forall \bm{d} \in \mathcal{M}_t\text{~and}$$
$$\max_{i\in \{1,\ldots,n_{\mathcal{X}}\}}\{\max_{s\in [t,T]}g_i(\bm{y}_{\bm{x},t}^{\bm{d}}(s),s)\}\leq 0, \forall \bm{d} \in \mathcal{M}_t.$$ Therefore, $$\bm{y}_{\bm{x},t}^{\bm{d}}(s) \in \mathcal{X}_s, \forall \bm{d}\in \mathcal{M}_t, \forall s \in [t,T]\text{~and}$$ $$\bm{y}_{\bm{x},t}^{\bm{d}}(T) \in \mathtt{TR},\forall \bm{d}\in \mathcal{M}_t.$$ Thus, $\bm{x} \in \mathcal{R}_t$.

Above all, $\mathcal{R}_t=\{\bm{x}\in \mathbb{R}^n\mid v(\bm{x},t)\leq 0\}$ and thus $\mathcal{R}_t=\{\bm{x}\in \mathbb{R}^n\mid u(\bm{x},t)\leq 0\}$. \hfill $\Box$
\end{IEEEproof}

According to Proposition \ref{relation}, the backward reachable set $\mathcal{R}_t$ is equal to the zero sub-level set of the value function $u(\bm{x},t)$ in \eqref{u}. In the following we show that this value function $u(\bm{x},t)$ is the unique Lipschitz continuous viscosity solution to the equation:
\begin{equation}
\label{HJB1}
\begin{split}
&\max\big\{\partial_t u(\bm{x},t)+H(\bm{x},\triangledown_{\bm{x}}u),\max_{i\in \{1,\ldots,n_{\mathcal{X}}\}} \{g_i(\bm{x},t)\}-u(\bm{x},t)\big\}=0
\end{split}
\end{equation}
with terminal condition
\[u(\bm{x},T)=\max\big\{\max_{i\in \{1,\ldots,n_{\mathtt{TR}}\}} \{l_i(\bm{x})\}, \max_{i\in \{1,\ldots,n_{\mathcal{X}}\}}\{g_{i}(\bm{x},T)\}\big\},\]
where $H(\bm{x},\bm{p})=\max_{\bm{d}\in \mathcal{D}}\bm{p}\cdot \bm{F}(\bm{x},\bm{d})$.

\eqref{HJB1} could be regarded as a special case of the Hamilton-Jacobi partial differential equation (4) in \cite{margellos2011}. \cite{margellos2011} considered reachability problems with competing inputs and time-invariant state constraints. In this paper we additionally consider time-varying state constraints. For the sake of clear presentation, in the following we give a brief introduction of inferring the HJE  \eqref{HJB1}. The viscosity solution $u(\bm{x},t)$ to \eqref{HJB1} is formalized in Definition \ref{viscosity}.
\begin{definition}\cite{Bardi1997,fialho1999}
\label{viscosity}
A lower semi-continuous function $u(\bm{x},t)$ on $\mathbb{R}^n\times [0,T]$ is called to be a viscosity super-solution of \eqref{HJB1}, if  for any test function $\psi\in C^{\infty}(\mathbb{R}^n\times [0,T])$ such that $u-\psi$ attains a local minimum at $(\bm{y}_0, t_0) \in \mathbb{R}^n\times [0,T]$,
\begin{equation}
\label{super}
\begin{split}
\max\big\{\partial_t &\psi(\bm{y}_0,t_0)+H(\bm{y}_0,\triangledown_{\bm{x}}\psi), \max_{i\in \{1,\ldots,n_{\mathcal{X}}\}}\{g_i(\bm{y}_0,t_0)\}-u(\bm{y}_0,t_0)\big\}\leq 0
\end{split}
\end{equation}
holds; A upper semi-continuous function $u(\bm{x},t)$ on $\mathbb{R}^n\times [0,T]$ is called to be a viscosity sub-solution of \eqref{HJB1}, if
for any test function $\psi\in C^{\infty}(\mathbb{R}^n\times [0,T])$ such that $u-\psi$ attains a local maximum at $(\bm{y}_0, t_0) \in \mathbb{R}^n\times [0,T]$,
\begin{equation}
\label{sub}
\begin{split}
\max\big\{\partial_t &\psi(\bm{y}_0,t_0)+H(\bm{y}_0,\triangledown_{\bm{x}}\psi), \max_{i\in \{1,\ldots,n_{\mathcal{X}}\}}\{g_i(\bm{y}_0,t_0)\}-u(\bm{y}_0,t_0)\big\}\geq 0
\end{split}
\end{equation}
holds.
A continuous function $u(\bm{x},t)$ on $\mathbb{R}^n\times [0,T]$ is called to be a viscosity solution to \eqref{HJB1} if it is both a viscosity super- and sub-solution to \eqref{HJB1}.
\end{definition}

Firstly, we show that the value function $u(\bm{x},t)$ is Lipschitz continuous, which is formally stated in Lemma \ref{p1}. 
    
\begin{lemma}
\label{p1}
Let $g_i$, $i=1,\ldots,n_{\mathcal{X}}$, and $l_j$, $j=1,\ldots,n_{\mathtt{TR}}$, be locally Lipschitz continuous functions respectively. Then $u(\bm{x},t)$ is locally Lipschitz continuous over $\mathbb{R}^n \times [0,T]$.
\end{lemma}
\begin{IEEEproof}
The proof is given in Appendix. \hfill $\Box$
\end{IEEEproof}

Secondly, $u(\bm{x},t)$ satisfies the dynamic programming principle presented in Lemma \ref{dp}.
\begin{lemma}
\label{dp}
For $(\bm{x},t)\in \mathbb{R}^n\times [0,T]$ and $\delta\geq 0$ satisfying $t+\delta\leq T$,
\begin{equation}
\label{dy}
\begin{split}
u(\bm{x},t)=&\sup_{\bm{d}\in \mathcal{M}_{[t,t+\delta]}}\max\big\{u(\bm{\phi}_{\bm{x},t}^{\bm{d}}(t+\delta),t+\delta),\max_{i\in \{1,\ldots,n_{\mathcal{X}}\}}\{\max_{s\in [t,t+\delta]}g_{i}(\bm{\phi}_{\bm{x},t}^{\bm{d}}(s),s)\}\big\},
\end{split}
\end{equation}
where $\bm{d}\in \mathcal{M}_{[t,t+\delta]}$ is the restriction of $\bm{d}\in \mathcal{M}_{t}$ over $[t,t+\delta]$.
\end{lemma}
\begin{IEEEproof}
The proof is given in Appendix. \hfill $\Box$
\end{IEEEproof}

We now show that the value function $u(\bm{x},t)$ in \eqref{u} is the unique continuous viscosity solution to \eqref{HJB1}.

\begin{theorem}
\label{solution}
The value function $u(\bm{x},t):\mathbb{R}^n\times [0,T]\mapsto \mathbb{R}$ in \eqref{u} is the unique Lipschitz continuous viscosity solution to HJE \eqref{HJB1}.
\end{theorem}
\begin{IEEEproof}
The proof is shown in Appendix. \hfill $\Box$
\end{IEEEproof}

We have shown that the value function $u(\bm{x},t)$ in \eqref{u}, whose zero sub-level set is the backward reachable set $\mathcal{R}_t$ at time $t\in [0,T]$, is
the unique Lipschitz continuous viscosity solution to HJE \eqref{HJB1}. Nowadays there are many efficient numerical methods for solving  \eqref{HJB1} with appropriate number of variables, e.g., \cite{bokanowski2013,falcone2016}. However, solving \eqref{HJB1} generally requires gridding the state space and thus is computationally intensive for some cases, especially for high-dimensional systems. In the section what follows, we will approximate this value function using polynomial viscosity super-solutions to \eqref{HJB1} by solving semi-definite programs. The Lipschitz continuity property of the viscosity solution to \eqref{HJB1} plays an important role in guaranteeing the existence of solutions to the constructed semi-definite program.

\section{Computing Inner Approximations}
\label{IACR}
In this section, by resorting to polynomial viscosity super-solutions to \eqref{HJB1} whose zero sub-level sets are inner-approximations of the backward reachable set $\mathcal{R}_0$, we first formulate the problem of computing inner approximations of the backward reachable set $\mathcal{R}_0$ as a semi-definite programming problem. We then prove that the interior of the backward reachable set $\mathcal{R}_0$ could be approximated in measure as the degree of the polynomial viscosity super-solutions tends to infinity under appropriate assumptions.

\subsection{Semi-definite Programming Implementation}
\label{SSP}
In this subsection we show that the zero sub-level set of a smooth viscosity super-solution to \eqref{HJB1} is an inner-approximation of the backward reachable set $\mathcal{R}_0$. Such a viscosity super-solution is computed by solving a semi-definite program, which is constructed from \eqref{HJB1}.

Firstly, we demonstrate that a smooth viscosity super-solution $\psi(\bm{x},t)$ to \eqref{HJB1} over $\mathbb{R}^n\times [0,T]$ is a solution to the following constraint:
\begin{equation}
\label{super2}
\begin{split}
\max\big\{\partial_t &\psi(\bm{x},t)+H(\bm{x},\triangledown_{\bm{x}}\psi), \max_{i\in \{1,\ldots,n_{\mathcal{X}}\}}\{g_i(\bm{x},t)\}-\psi(\bm{x},t)\big\}\leq 0.
\end{split}
\end{equation}
This conclusion is stated in Lemma \ref{con_super} formally.
\begin{lemma}
\label{con_super}
Assume that $\psi\in C^{\infty}(\mathbb{R}^n\times [0,T])$. $\psi$ is a viscosity super-solution of \eqref{HJB1} if and only if it satisfies the constraint \eqref{super2} over $\mathbb{R}^n\times [0,T]$.
\end{lemma}
\begin{IEEEproof}
First, we prove that if $\psi(\bm{x},t)$ is a viscosity super-solution of \eqref{HJB1}, it satisfies \eqref{super2}.

According to the definition of the viscosity super-solution in Definition \ref{viscosity}, i.e. for all test function $v\in C^{\infty}(\mathbb{R}^n\times [0,T])$ such that $\psi-v$ attains a local minimum  at $(\bm{x}_0,t_0)$,
then
\begin{equation}
\begin{split}
\max\{\partial_s &v(\bm{x}_0,t_0)+H(\bm{x}_0,\triangledown_{\bm{x}}v),\max_{i\in \{1,\ldots,n_{\mathcal{X}}\}}\{g_i(\bm{x}_0,t_0)\}-\psi(\bm{x}_0,t_0)\}\leq 0.
\end{split}
\end{equation}
It is apparent that $\psi \in C^{\infty}(\mathbb{R}^n\times [0,T])$ satisfies \eqref{super2} since $\psi-v=\psi-\psi$ when $v=\psi$ attains local minimum at any $(\bm{x}_0,t_0) \in \mathbb{R}^n\times [0,T]$.

Next, we prove that if $\psi(\bm{x},t)$ satisfies \eqref{super2}, it is a viscosity super-solution of \eqref{HJB1}. This claim can be assured by following the proof of Theorem \ref{solution} for the viscosity super-solution part. Let $v\in C^{\infty}(\mathbb{R}^n\times [0,T])$ such that $\psi-v$ attains a local minimum at $(\bm{x}_0,t_0)$, where $t_0\in [0,T]$. Similarly, we assume that this minimum is $0$, i.e. $v(\bm{x}_0,t_0)=\psi(\bm{x}_0,t_0)$.

If \eqref{super} is false, then either
\begin{equation}
\label{con211}
\max_{i\in \{1,\ldots,n_{\mathcal{X}}\}}\{g_i(\bm{x}_0,t_0)\}\geq v(\bm{x}_0,t_0)+\epsilon_1
\end{equation}
holds
or
\begin{equation}
\label{con221}
\partial_t v(\bm{x}_0,t_0)+H(\bm{x}_0,\nabla_{\bm{x}}v)\geq \epsilon_2
\end{equation}
holds for some $\epsilon_1$, $\epsilon_2>0$.

If \eqref{con211} holds, then there is a small enough $\delta>0$ such that for $(\bm{x},t)$ satisfying $t_0\leq t\leq t_0+\delta$ and
$\|\bm{x}-\bm{x}_0\|\leq \delta$,
\[\max_{i\in \{1,\ldots,n_{\mathcal{X}}\}}\{g_i(\bm{x},t)\}\geq v(\bm{x}_0,t_0)+\frac{\epsilon_1}{2}=\psi(\bm{x}_0,t_0)+\frac{\epsilon_1}{2}.\]
However, since $\psi$ satisfies \eqref{super2} over $\mathbb{R}^n\times [0,T]$, we have $\psi(\bm{x},t)\geq \max_{i\in \{1,\ldots,n_{\mathcal{X}}\}}\{g_i(\bm{x},t)\}$, implying that $\psi(\bm{x}_0,t_0)\geq \psi(\bm{x}_0,t_0)+\frac{\epsilon_1}{2}$,
which is a contradiction since $\epsilon_1>0$.

However, if \eqref{con221} holds, there is a small enough $\delta>0$ such that there exists a strategy $\bm{d}_1\in \mathcal{M}_t$ such that
\begin{equation}
\label{exist1}
\delta\frac{\epsilon_2}{2}\leq v(\bm{\phi}_{\bm{x}_0,t_0}^{\bm{d}_1}(t_0+\delta),t_0+\delta)-v(\bm{x}_0,t_0).
\end{equation}
Further, due to the fact that 
\begin{equation}
\delta\frac{\epsilon_2}{2}\leq \psi(\bm{\phi}_{\bm{x}_0,t_0}^{\bm{d}_1}(t_0+\delta),t_0+\delta)-\psi(\bm{x}_0,t_0),
\end{equation}
which contradicts $\psi(\bm{\phi}_{\bm{x}_0,t_0}^{\bm{d}_1}(t_0+\delta),t_0+\delta)\leq \psi(\bm{x}_0,t_0)$, which is obtained by the fact that $\partial_s \psi(\bm{x}_0,t_0)+H(\bm{x}_0,\triangledown_{\bm{x}}\psi)\leq 0$ and $\psi\in C^{\infty}(\mathbb{R}^n\times [0,T])$. Therefore, if $\psi(\bm{x},t)$ satisfies \eqref{super2}, it is a viscosity super-solution of \eqref{HJB1}.

Therefore, the conclusion holds. \hfill $\Box$
 \end{IEEEproof}

Based on Lemma \ref{con_super} we will show that an inner-approximation of the backward reachable set $\mathcal{R}_0$ can be characterized by the zero sub-level set of a smooth viscosity super-solution to \eqref{HJB1}. 
\begin{theorem}
\label{inner}
If $\psi(\bm{x},t)\in C^{\infty}(\mathbb{R}^n\times [0,T])$ is a viscosity super-solution of \eqref{HJB1} with boundary condition $$\psi(\bm{x},T)\geq \max\big\{\max_{i\in \{1,\ldots,n_{\mathtt{TR}}\}}\{l_i(\bm{x})\},\max_{i\in \{1,\ldots,n_{\mathcal{X}}\}}\{g_i(\bm{x},T)\}\big\}$$ over $\bm{x}\in \mathbb{R}^n,$ $\{\bm{x}\in \mathbb{R}^n\mid \psi(\bm{x},0)\leq 0\}$ is an inner-approximation of the backward reachable set $\mathcal{R}_0$.
\end{theorem}
\begin{IEEEproof}
For each $\bm{d}\in \mathcal{D}$ and $\bm{x}_0\in \mathbb{R}^n$,
 $$\psi(\bm{\phi}_{\bm{x}_0,0}^{\bm{d}}(t),t)-\psi(\bm{x}_0,0)=\int_{0}^t\mathcal{L}_{\bm{d}}\psi(\bm{x},s)ds$$
 for $t\in [0,T]$,
where $\mathcal{L}_{\bm{d}}\psi(\bm{x},s)=\partial_s \psi(\bm{x},s)+\nabla_{\bm{x}}\psi \cdot \bm{F}(\bm{x},\bm{d}).$ According to Lemma \ref{con_super}, $\mathcal{L}_{\bm{d}}\psi(\bm{x},s)\leq 0$ for $s\in [0,t]$ holds. Therefore,
\[\psi(\bm{\phi}_{\bm{x}_0,0}^{\bm{d}}(t),t)\leq \psi(\bm{x}_0,0) \] for $t\in [0,T]$.
Obviously, if $\psi(\bm{x}_0,0)\leq 0$, then $\psi(\bm{\phi}_{\bm{x}_0,0}^{\bm{d}}(T),T)\leq 0$ holds. Since $$\psi(\bm{x},T)\geq \max\big\{\max_{i\in \{1,\ldots,n_{\mathtt{TR}}\}}\{l_i(\bm{x})\},\max_{i\in \{1,\ldots,n_{\mathcal{X}}\}}\{g_i(\bm{x},T)\}\big\}$$ over $\bm{x}\in \mathbb{R}^n$ and $\{\bm{x}\in \mathbb{R}^n\mid \max\big\{\max_{i\in \{1,\ldots,n_{\mathtt{TR}}\}}\{l_i(\bm{x})\},\max_{i\in \{1,\ldots,n_{\mathcal{X}}\}}\{g_i(\bm{x},T)\}\big\}\leq 0\}\subset \mathtt{TR}$, $$\bm{\phi}_{\bm{x}_0,0}^{\bm{d}}(T) \in \mathtt{TR}$$ holds. Also, according to Lemma \ref{con_super}, $\psi(\bm{x},t)\geq \max_{i\in \{1,\ldots,n_{\mathcal{X}}\}}\{g_i(\bm{x},t)\}$ over $(\bm{x},t)\in \mathbb{R}^n\times [0,T]$. Since $\mathcal{X}_t=\{\bm{x}\in \mathbb{R}^n\mid \max_{i\in \{1,\ldots,n_{\mathcal{X}}\}} \{g_i(\bm{x},t)\}\leq 0\}$ for $t\in [0,T]$, if $\psi(\bm{x}_0,0)\leq 0$, then we have $\bm{\phi}_{\bm{x}_0,0}^{\bm{d}}(t)\in \mathcal{X}_t$ for $ t\in [0,T]$. Therefore, $\psi(\bm{x}_0,0)\leq 0$ implies that  all trajectories starting from $\bm{x}_0$ will enter the target region $\mathtt{TR}=\{\bm{x}\in \mathbb{R}^n\mid \max_{i\in \{1,\ldots,n_{\mathtt{TR}}\}}\{l_i(\bm{x})\}\leq 0\}$ after the time duration $T$ while staying inside the constraint set  $\mathcal{X}_t$ over $t\in [0,T]$. Therefore, $\{\bm{x}\in \mathbb{R}^n\mid \psi(\bm{x},0)\leq 0\}$ is an inner-approximation of the backward reachable set $\mathcal{R}_0$. \hfill $\Box$
\end{IEEEproof}

According to Theorem \ref{inner} and Lemma \ref{con_super}, the problem of computing a smooth viscosity super-solution $\psi(\bm{x},s)$, whose zero sub-level set is an inner-approximation of the backward reachable set $\mathcal{R}_0$, can be reformulated as the following constraint over $\psi\in C^{\infty}(\mathbb{R}^n\times [0,T])$,
\begin{equation}
\label{optimization}
\begin{split}
&\partial_t \psi(\bm{x},t)+\triangledown_{\bm{x}}\psi \cdot \bm{F}(\bm{x},\bm{d})\leq 0, \forall (t,\bm{x},\bm{d})\in [0,T]\times \mathbb{R}^n\times \mathcal{D},\\
&\psi(\bm{x},t)-g_i(\bm{x},t)\geq 0,~~~~~~~~~~~\forall (\bm{x},s)\in \mathbb{R}^n\times [0,T],\\
&\psi(\bm{x},T)-l_j(\bm{x})\geq 0,~~~~~~~~~~~~~\forall \bm{x} \in \mathbb{R}^n,\\
&i=1,\ldots,n_{\mathcal{X}}, j=1,\ldots,n_{\mathtt{TR}}.\\
\end{split}
\end{equation}

\begin{corollary}
\label{inner-op}
Let $\psi(\bm{x},t)$ be a solution to \eqref{optimization}. Then $\{\bm{x}\in \mathbb{R}^n\mid \psi(\bm{x},0)\leq 0\}$ is an inner-approximation of the backward reachable set $\mathcal{R}_0$.
\end{corollary}
\begin{IEEEproof}
The claim in this corollary can be easily assured by Lemma \ref{con_super} and Theorem \ref{inner}. \hfill $\Box$
\end{IEEEproof}

The problem of obtaining a solution to \eqref{optimization} is challenging since a solution $\psi(\bm{x},t)$ should satisfy \eqref{optimization} for $\bm{x}\in \mathbb{R}^n$. In the following we will relax this condition to obtain a function $\psi(\bm{x},t)$ satisfying \eqref{optimization} for $\bm{x}\in B(\bm{0},R)$, where $B(\bm{0},R)$ is defined in \eqref{B}.

Regarding $\bm{x}\in B(\bm{0},R)$ and $\bm{d}\in \mathcal{D}$, we have $\bm{F}(\bm{x},\bm{d})=\bm{f}(\bm{x},\bm{d})$. When the viscosity super-solution $\psi(\bm{x},t)$ to \eqref{optimization} is constrained to polynomial type in the set $B(\bm{0},R)\times [0,T]$, \eqref{optimization} can be recast as the sum-of-squares program \eqref{sos}, which is formalized below. The constraints in \eqref{sos} that polynomials are sum-of-squares can be written explicitly as LMI constraints, and the objective is linear in the coefficients of the polynomial $\psi(\bm{x},t)$; therefore problem \eqref{sos} is able to be formulated as a semi-definite program, which falls within the convex programming framework and can be solved via interior-point methods in polynomial time. Note that the objective of \eqref{sos} would facilitate the gain of less conservative inner-approximations of the backward reachable set. The reason is that if $\psi_1(\bm{x},0)\leq \psi_2(\bm{x},0)$ over $\bm{x}\in B(\bm{0},R)$, then $\{\bm{x}\in B(\bm{0},R)\mid \psi_2(\bm{x},0)\leq 0\}\subseteq \{\bm{x}\in B(\bm{0},R)\mid \psi_1(\bm{x},0)\leq 0\}$ and $\int_{B(\bm{0},R)} \psi_1(\bm{x},0)d\bm{x}\leq \int_{B(\bm{0},R)} \psi_2(\bm{x},0)d\bm{x}$.
\begin{algorithm}
\begin{algorithmic}
\STATE
\begin{equation}
\label{sos}
\begin{split}
&d_k^*=\inf~\bm{w}'\cdot \bm{l} \\
&s.t.\\
&-\mathcal{L}\psi(\bm{x},t)=s_0+s_1g_R(\bm{x})+s_2t(T-t)+\sum_{i=1}^{n_{\mathcal{D}}}s'_{i}h_i(\bm{d}),\\
&\psi(\bm{x},t)-g_i(\bm{x},t)=s_{3,i}+s_{4,i}g_R(\bm{x})+s_{5,i}t(T-t),\\
&\psi(\bm{x},T)-l_j(\bm{x})=s_{6,j}+s_{7,j}g_R(\bm{x}),\\
&i=1,\ldots,n_{\mathcal{X}},j=1,\ldots,n_{\mathtt{TR}},\\
\end{split}
\end{equation}
where $\bm{w}'\cdot \bm{l}=\int_{B(\bm{0},R)} \psi(\bm{x},0)d\bm{x}$, $\bm{l}$ is the vector of the moments of the Lebesgue measure over $B(\bm{0},R)$ indexed in the same basis in which the polynomial $\psi(\bm{x},0)$ with coefficients $\bm{w}$ is expressed, and $\mathcal{L}\psi(\bm{x},t)=\partial_t\psi(\bm{x},t)+\nabla_{\bm{x}}\psi(\bm{x},t)\cdot\bm{f}(\bm{x},\bm{d})$. The minimum is over polynomial $\psi(\bm{x},t) \in \mathbb{R}_{k}[\bm{x},t]$ and sum-of-squares polynomials $s_0(\bm{x},t,\bm{d})$, $s_1(\bm{x},t,\bm{d})$, $s_2(\bm{x},t,\bm{d})$, $s'_r(\bm{x},t,\bm{d})$, $r=1,\ldots,n_{\mathcal{D}}$, $s_{3,i}(\bm{x},t)$, $s_{4,i}(\bm{x},t)$, $s_{5,i}(\bm{x},t)$, $i=1,\ldots,n_{\mathcal{X}}$ and $s_{6,j}(\bm{x})$, $s_{7,j}(\bm{x})$, $j=1,\ldots, n_{\mathtt{TR}}$, of appropriate degrees.
\end{algorithmic}
\label{alg}
\end{algorithm}

\begin{corollary}
\label{inner2}
Let $\psi(\bm{x},t)$ be a solution to \eqref{sos}, then $\{\bm{x}\in B(\bm{0},R)\mid \psi(\bm{x},0)\leq 0\}$ is an inner-approximation of the backward reachable set $\mathcal{R}_0$.
\end{corollary}
\begin{IEEEproof}
According to \eqref{sos}, we have the following constraints: for
$(\bm{x},t,\bm{d})\in B(\bm{0},R)\times [0,T] \times \mathcal{D}$,
\begin{equation}
\label{optimization1}
\begin{split}
&\partial_t \psi(\bm{x},t)+\triangledown_{\bm{x}}\psi \cdot \bm{f}(\bm{x},\bm{d})\leq 0,\\
&\psi(\bm{x},t)-g_i(\bm{x},t)\geq 0,\\
&\psi(\bm{x},T)-l_j(\bm{x})\geq 0,\\
&i=1,\ldots,n_{\mathcal{X}}, j=1,\ldots,n_{\mathtt{TR}}.\\
\end{split}
\end{equation}
Obviously, $\{\bm{x}\in B(\bm{0},R)\mid \psi(\bm{x},0)\leq 0\}\subseteq \mathcal{X}_0$.
Assume that there exists a trajectory initialized in $\bm{x}_0\in \{\bm{x}\in B(\bm{0},R)\mid \psi(\bm{x},0)\leq 0\}$ at $t=0$ such that it escapes from the set $\mathcal{X}_{\tau}$ at some $\tau\in [0,T]$, i.e. $\bm{\phi}_{\bm{x}_0,0}^{\bm{d}_1}(\tau) \notin \mathcal{X}_{\tau}$ for some $\bm{d}_1\in \mathcal{D}$. Therefore, there exists $\tau_1\in [0,\tau]$ such that $\bm{\phi}_{\bm{x}_0,0}^{\bm{d}_1}(\tau_1)\notin \mathcal{X}_{\tau_1}$ and $\bm{\phi}_{\bm{x}_0,0}^{\bm{d}_1}(\tau_1)\in B(\bm{0},R)$, implying that $\psi(\bm{\phi}_{\bm{x}_0,0}^{\bm{d}_1}(\tau_1),\tau_1)>0$. This contradicts the  fact that  $ \psi(\bm{x}_0,0)\leq 0$ and $\partial_t \psi(\bm{x},t)+\triangledown_{\bm{x}}\psi \cdot \bm{f}(\bm{x},\bm{d})\leq 0$. We conclude that every possible trajectory originating in $\{\bm{x}\in B(\bm{0},R)\mid \psi(\bm{x},0)\leq 0\}$ at $t=0$ will stay inside the set $\mathcal{X}_t$, $\forall t\in [0,T]$. Also, since
$\psi(\bm{x},T)-l_j(\bm{x})\geq 0, j=1,\ldots,n_{\mathtt{TR}}$, by following the proof of Theorem \ref{inner}, we conclude that $\{\bm{x}\in B(\bm{0},R)\mid \psi(\bm{x},0)\leq 0\}\subseteq \mathcal{R}_0$. \hfill $\Box$
\end{IEEEproof}

Corollary \ref{inner2} implies that an inner-approximation of the backward reachable set $\mathcal{R}_0$ is able to be synthesized by solving the semi-definite program \eqref{sos}. In the following subsection we prove the existence of a convergent sequence of inner-approximations, which are formed by solutions to \eqref{sos}, to the interior of the backward reachable set $\mathcal{R}_0$ in measure under appropriate assumptions.

\subsection{Convergence Analysis}
\label{CA}
In this subsection we show that \eqref{sos} exhibits a convergent sequence of inner approximations to the interior of the backward reachable set $\mathcal{R}_0$ in measure under appropriate assumptions. We firstly show that on a given compact set $B(\bm{0},R)$ there is a smooth solution to \eqref{optimization1}, which can approximate the viscosity solution $u$ to \eqref{HJB1} uniformly. Then we demonstrate that there exists a sequence of polynomial functions satisfying \eqref{sos} and approximating the viscosity solution $u$ uniformly.

Before this, we introduce an auxiliary lemma stating that over a compact set $B(\bm{0},R)\times [0,T]$ there is a smooth function $\psi(\bm{x},t)$ which is an uniform approximation to a given Lipschitz function $u(\bm{x},t)$ and provides one side approximation to the Dini$-$derivative of the form
$\sup_{\bm{d}\in \mathcal{D}}\mathcal{L} \psi(\bm{x},t)$,
where $\mathcal{L}\psi(\bm{x},t)=\partial_t \psi+\nabla_{\bm{x}} \psi\cdot \bm{f}(\bm{x},\bm{d})$.

\begin{lemma}\cite{lin1996}
\label{smooth}
Let $B(\bm{0},R)$ be a compact subset in $\mathbb{R}^n$ and $u(\bm{x},t):  B(\bm{0},R)\times [0,T] \mapsto \mathbb{R}$ be a Lipschitz function. If there exists a continuous function $\alpha: B(\bm{0},R)\times [0,T] \mapsto \mathbb{R}$ such that for each $\bm{d}\in \mathcal{D}$,
\[\mathcal{L} u(\bm{x},t)\leq \alpha(\bm{x},t), \text{a.e.} (\bm{x},t) \in B(\bm{0},R)\times [0,T],\]
(recall that $\nabla_{\bm{x}} u$ is defined a.e., since $u$ is locally Lipschitz.)
 then for any given $\epsilon>0$, there exists some smooth function $\psi(\bm{x},t)$ defined on $B(\bm{0},R)\times [0,T]$ such that
\[\sup_{(\bm{x},t)\in B(\bm{0},R)\times [0,T]}|\psi(\bm{x},t)-u(\bm{x},t)|<\epsilon \text{~and}\]
\[\sup_{\bm{d}\in  \mathcal{D}}\mathcal{L} \psi(\bm{x},t)\leq \alpha(\bm{x},t)+\epsilon\]
over $(\bm{x},t)\in B(\bm{0},R)\times [0,T]$.
\end{lemma}

According to Lemma \ref{smooth}, we have the following theorem stating that there exists a smooth viscosity super-solution $\psi(\bm{x},t)$ with $\psi(\bm{x},t)\geq u(\bm{x},t)$ approximating the viscosity solution $u(\bm{x},t)$ to \eqref{HJB1} uniformly on the compact set $ B(\bm{0},R)\times [0,T]$.
\begin{theorem}
\label{smooth-super}
There exists a smooth function $\psi(\bm{x},t)$ satisfying \eqref{optimization1}, which approximates the viscosity solution $u(\bm{x},t)$ to \eqref{HJB1} uniformly on the compact set $B(\bm{0},R)\times [0,T]$.
\end{theorem}
\begin{IEEEproof}
From Theorem \ref{solution}, we have that the value function $u(\bm{x},t)$ is locally Lipschitz continuous. Also, since $\mathcal{L}u \leq 0$ for all $(\bm{x},t,\bm{d})\in B(\bm{0},R)\times [0,T]\times \mathcal{D}$, where $\mathcal{L}u(\bm{x},t)=\partial_t u+\nabla_{\bm{x}} u\cdot \bm{f}(\bm{x},\bm{d})$, according to Lemma \ref{smooth}, there exists a smooth function $\overline{\psi}: B(\bm{0},R)\times [0,T] \mapsto \mathbb{R}$ such that
\[\sup_{(t,\bm{x})\in [0,T]\times B(\bm{0},R)}|\overline{\psi}(\bm{x},t)-u(\bm{x},t)|<\epsilon\text{~and}\]
\[\sup_{(\bm{x},t,\bm{d})\in B(\bm{0},R)\times [0,T]\times  \mathcal{D}}\mathcal{L} \overline{\psi}(\bm{x},t)\leq \epsilon.\]
Let
\[\psi(\bm{x},t)= \overline{\psi}(\bm{x},t)-\epsilon t+(T+1)\epsilon.\]
Since $$u(\bm{x},t)-\epsilon\leq \overline{\psi}(\bm{x},t)\leq u(\bm{x},t)+\epsilon$$ for $ (\bm{x},t)\in B(\bm{0},R)\times [0,T]$,
$$u+\epsilon(T-t)\leq \psi \leq u+(T-t+2)\epsilon$$ holds for $ (\bm{x},t)\in B(\bm{0},R)\times [0,T]$. Thus $u\leq \psi$ for $(\bm{x},t)\in B(\bm{0},R)\times [0,T]$. Also, $\psi: B(\bm{0},R)\times [0,T]\mapsto \mathbb{R}$ approximates $u(\bm{x},t)$ uniformly on the compact set $B(\bm{0},R)\times [0,T]$.

We also need to prove that $\psi(\bm{x},t)\geq \max_{i\in \{1,\ldots,n_{\mathcal{X}}\}}\{g_i(\bm{x},t)\}$ and $\mathcal{L}\psi\leq 0$ over $(\bm{x},t,\bm{d})\in B(\bm{0},R)\times [0,T]\times \mathcal{D}$. The former, i.e. $\max_{i\in \{1,\ldots,n_{\mathcal{X}}\}}\{g(\bm{x},t)\}\leq \psi(\bm{x},t)$ over $(\bm{x},t)\in B(\bm{0},R)\times [0,T]$, holds obviously since $\max_{i\in \{1,\ldots,n_{\mathcal{X}}\}}\{g_i(\bm{x},t)\}\leq u$. The latter, i.e. $$\mathcal{L} \psi\leq 0$$ also holds since $$\mathcal{L}\psi\leq \mathcal{L} \overline{\psi}-\epsilon$$ for $ (\bm{x},t,\bm{d})\in  B(\bm{0},R)\times [0,T] \times \mathcal{D}$.

Thus, the conclusion in Theorem \ref{smooth-super} holds. \hfill $\Box$
\end{IEEEproof}

In our implementation we restrcit smooth viscosity super-solutions to polynomial functions and attempt to inner-approximate the backward reachable set $\mathcal{R}_0$ by solving the semi-definite programming problem \eqref{sos}. In the following we prove that under Assumption \ref{DD} there exists a sequence of polynomials $\{\psi_{k+N}(\bm{x},t)\}_{k\in \mathbb{N}}$, where $N$ is some positive integer and $\psi_{k+N}(\bm{x},t)\in \mathbb{R}_{k+N}[\bm{x},t]$ satisfies \eqref{sos},  such that $\lim_{k\rightarrow +\infty}\psi_{k+N}(\bm{x},t)=u(\bm{x},t)$  uniformly over $B(\bm{0},R)\times [0,T]$, where $u(\bm{x},t)$ is the unique viscosity solution to \eqref{HJB1}.

\begin{assumption}
\label{DD}
One of the polynomials defining the set $\mathcal{D}$ is equal to $R_{\mathcal{D}}-\|\bm{d}\|^2$ for some constant $R_{\mathcal{D}}\geq 0$.
\end{assumption}

Assumption \ref{DD} is without loss of generality because of compactness of $\mathcal{D}$. Thus $R_{\mathcal{D}}-\|\bm{d}\|^2\geq 0$ can be  a redundant constraint defining $\mathcal{D}$ for sufficiently large $R_{\mathcal{D}}$.

The proof also requires Putinar's Positivstellensatz.
\begin{theorem}
\label{putinar}
[Putinar's Positivstellensatz \cite{putinar93}] Let $\mathcal{K}=\{\bm{y}\in \mathbb{R}^m| g_1(\bm{y})\geq 0,\ldots,g_l(\bm{y})\geq 0\}$ be a compact set. Suppose there exists $N>0$ such that
\[N-\sum_{i=1}^m y_i^2\in M(g_1,\ldots,g_l).\]
If $p(\bm{y})$ is positive on $K$, then $p(\bm{y})\in M(g_1,\ldots,g_l),$ where $M(g_1,\ldots,g_l)$ is the quadratic module of polynomials $g_1,\ldots,g_l,$ i.e.
\[M(g_1,\ldots,g_l)=\{\sigma_0(\bm{y})+\sum_{i=1}^l \sigma_i(\bm{y})g_i(\bm{y})|~\text{each}~\sigma_i\in \sum_m\}\] with
$\sum_m$ being the set of sum of squares (SOS) polynomials over variables $\bm{y}$, i.e.
\[\sum_m:=\{p\in \mathbb{R}[\bm{y}]|p=\sum_{i=1}^k q_i^2,q_i\in \mathbb{R}[\bm{y}],i=1,\ldots,k\}.\]
\end{theorem}

\begin{corollary}
\label{unifrom}
There exists a sequence of polynomials $\{\psi_{k+N}(\bm{x},t)\}_{k=1}^{\infty}$ satisfying \eqref{sos}, where $\psi_{k+N}(\bm{x},t)\in \mathbb{R}_{k+N}[\bm{x},t]$ and $N$ is some positive integer, such that $\psi_{k+N}(\bm{x},t)$ approximates $u(\bm{x},t)$ over the space $B(\bm{0},R)\times [0,T]$ uniformly as $k$ approaches infinity, where $u(\bm{x},t)$ is the viscosity solution to \eqref{HJB1}. 
\end{corollary}
\begin{IEEEproof}
From Theorem \ref{smooth-super}, for every $\epsilon>0$, there exists a smooth function  $v(\bm{x},t)\in C^{\infty}(B(\bm{0},R)\times [0,T])$ such that $$v(\bm{x},t)\geq u(\bm{x},t), $$ $$|v(\bm{x},t)-u(\bm{x},t)|\leq \epsilon, \text{~and}$$
$$\mathcal{L}v\leq 0$$
for $ (\bm{x},t,\bm{d})\in B(\bm{0},R)\times [0,T]\times \mathcal{D}$. 
Next, we follow the proof of Theorem 5 in \cite{henrion2014}.
Let
\[\tilde{v}(\bm{x},t)=v(\bm{x},t)-\epsilon t+(T+1)\epsilon.\]
We have $\mathcal{L}\tilde{v}=\mathcal{L} v-\epsilon\leq -\epsilon$, $\max_{i\in \{1,\ldots,n_{\mathcal{X}}\}}\{g_i(\bm{x},t)\}-\tilde{v}(\bm{x},t)\leq -\epsilon$ and $\tilde{v}(\bm{x},T)-\max_{i\in \{1,\ldots,n_{\mathtt{TR}}\}}\{l_i(\bm{x})\}=v(\bm{x},T)+\epsilon-\max_{i\in \{1,\ldots,n_{\mathtt{TR}}\}}\{l_i(\bm{x})\}\geq \epsilon$.
Since $[0,T]\times B(\bm{0},R)$ is compact, there exists a polynomial $\psi_N \in \mathbb{R}_N[\bm{x},t]$ of a sufficiently high degree $N$ such that
\[\sup_{[0,T]\times B(\bm{0},R)}|\tilde{v}-\psi_N|<\epsilon \text{~and}\]
\[\sup_{[0,T]\times B(\bm{0},R)\times \mathcal{D}} |\mathcal{L}\tilde{v}-\mathcal{L}\psi_N|<\epsilon.\]  The polynomial $\psi_N$ is therefore strictly feasible in \eqref{sos} (this follows from the classical Putinar's Positivstellensatz, as formulated in Theorem \ref{putinar}), and moreover $\psi_N(\bm{x},t)\geq u(\bm{x},t)$ for  $(\bm{x},t) \in  B(\bm{0},R)\times [0,T]$. 
Also,
$$\sup_{[0,T]\times B(\bm{0},R)}|\psi_N-u|<\epsilon (T+3).$$

Since $\epsilon$ is arbitrary, we conclude that as the degree $k$ tends to infinity, $\psi_{k+N}(\bm{x},t)\in \mathbb{R}_{k+N}[\bm{x},t]$ converges to $u(\bm{x},t)$ uniformly over $B(\bm{0},R)\times [0,T]$. \hfill $\Box$
\end{IEEEproof}

Corollary \ref{unifrom} establishes a uniformly functional convergence of $\psi$ to $u$. Let the sequence of polynomials $\{\psi_{k+N}(\bm{x},t)\}_{k=0}^{+\infty}$ satisfy Corollary \ref{unifrom}, and
\[\mathcal{X}_{0k}=\{\bm{x}\in \mathcal{X}_0\mid \psi_{k+N}(\bm{x},0)\leq 0\}.\]
Obviously, the inclusion $\mathcal{X}_{0k}\subset \mathcal{R}_0$ holds for all $k\in \{1,2,\ldots\}.$ Finally, we will show that $\mathcal{X}_{0k}$ approximates the interior of the backward reachable set $\mathcal{R}_0$ in measure. Before this, we first prove the non-fattening property of the zero level set of $u(\bm{x},t)$ evolving over time.

\begin{lemma}
\label{zero}
$\{\bm{x}\in \mathbb{R}^n\mid  u(\bm{x},0)=0\}$ is the boundary of the set
$\{\bm{x}\in \mathbb{R}^n\mid u(\bm{x},0)\leq 0\}$, where $u(\bm{x},t)$ is the viscosity solution to \eqref{HJB1}.
\end{lemma}
\begin{IEEEproof}
The fact that $\{\bm{x}\in \mathbb{R}^n\mid u(\bm{x},T)=0\}=\partial \{\bm{x}\in \mathbb{R}^n\mid u(\bm{x},T)\leq 0\}$ is easily assured by the fact that $u(\bm{x},T)=\max\{\max_{i\in \{1,\ldots,n_{\mathtt{TR}}\}}l_i(\bm{x}), \max_{i\in \{1,\ldots,n_{\mathcal{X}}\}}g_i(\bm{x},T)\}$, $\partial \mathcal{X}_t=\cup_{i=1}^{n_{\mathcal{X}}} \{\bm{x}\in \mathcal{X}_t\mid g_i(\bm{x},t)=0\}$ for $t\in [0,T]$ and $\partial \mathtt{TR}=\cup_{i=1}^{n_{\mathtt{TR}}} \{\bm{x}\in \mathtt{TR}\mid l_i(\bm{x})=0\}$.

Suppose that $\bm{y}\in \{\bm{x} \in \mathbb{R}^n \mid u(\bm{x},T)=0\}$ and there exists $\bm{y}_0$ belonging to the interior of the set $\{\bm{x} \in \mathbb{R}^n \mid u(\bm{x},0)\leq 0\}$ and $\bm{d}\in \mathcal{M}_0$ such that $\bm{y}=\bm{\phi}_{\bm{y}_0,0}^{\bm{d}}(T)$. Thus, there exist $\bm{y}_{1}$ and $\bm{x}_1$ satisfying $\bm{y}_{1}\notin \{\bm{x}\in \mathbb{R}^n\mid u(\bm{x},T)\leq 0\}$ but $\bm{y}_1\in \mathbb{U}(\bm{y};\epsilon)$ and $\bm{x}_1\in \mathbb{U}(\bm{y}_0;\epsilon_1) \subset \{\bm{x}\in \mathbb{R}^n\mid u(\bm{x},0)\leq 0\}$ such that $\bm{y}_1=\bm{\phi}_{\bm{x}_1,0}^{\bm{d}}$, where $\mathbb{U}(\cdot;\delta)$ with $\delta>0$ denotes the $\delta-$neighborhood of the argument. Thus this contradicts the fact that the states $\bm{x}$s in $\{\bm{x} \in \mathbb{R}^n\mid u(\bm{x},0)\leq 0\}$ will enter the target set $\{\bm{x}\in \mathbb{R}^n\mid u(\bm{x},T)\leq 0\}$ after the time duration of $T$ for $\bm{d} \in \mathcal{M}_0$. Thus, $\bm{y}_0$ belongs to the boundary of the set $\{\bm{x} \in \mathbb{R}^n \mid u(\bm{x},0)\leq 0\}$, implying $u(\bm{y}_0,0)=0$.

Since $\partial \{\bm{x}\in \mathbb{R}^n \mid u(\bm{x},0)\leq 0\}\subseteq \{\bm{x}\in \mathbb{R}^n\mid u(\bm{x},0)=0\}$ is clear, it is sufficient to prove that $\{\bm{x}\in \mathbb{R}^n\mid u(\bm{x},0)=0\}\subseteq \partial \{\bm{x}\in \mathbb{R}^n \mid u(\bm{x},0)\leq 0\}$. Let $u(\bm{x}_0,0)=0$. Since pointwise limits of measurable functions are measurable, $\mathcal{M}_0$ is a closed subset and thus remains compact \cite{wal95,platzer2017}. Therefore, there exists $\bm{d}_{\bm{x}_0}\in \mathcal{M}_0$ such that
\begin{equation}
\begin{split}
u(\bm{x}_0,0)=&\max\Big\{\max_{i\in \{1,\ldots,n_{\mathtt{TR}}\}}\{l_i(\bm{\phi}_{\bm{x}_0,0}^{\bm{d}_{\bm{x}_0}}(T))\},\max_{i\in \{1,\ldots,n_{\mathcal{X}}\}}\{\max_{s\in [0,T]}g_i(\bm{\phi}_{\bm{x}_0,0}^{\bm{d}_{\bm{x}_0}}(s),s)\}\Big\}.
\end{split}
\end{equation}
 We will prove that $\bm{x}_0\in \partial \{\bm{x}\in \mathbb{R}^n\mid u(\bm{x},0)\leq 0\}.$ 

 Assume that $\bm{x}_0$ belongs to the interior of the set $\{\bm{x}\in \mathbb{R}^n\mid u(\bm{x},0)\leq 0\}$. Then $$\max_{i\in \{1,\ldots,n_{\mathcal{X}}\}}\{g_i(\bm{\phi}_{\bm{x}_0,0}^{\bm{d}_{\bm{x}_0}}(s),s)\}<0.$$ Moreover, from above discussions, we deduce that
\begin{equation}
\begin{split}
&\max \big\{\max_{i\in \{1,\ldots,n_{\mathtt{TR}}\}}\{l_i(\bm{\phi}_{\bm{x}_0,0}^{\bm{d}_{\bm{x}_0}}(T))\},\max_{i\in \{1,\ldots,n_{\mathcal{X}}\}}\{g_i(\bm{\phi}_{\bm{x}_0,0}^{\bm{d}_{\bm{x}_0}}(T),T)\}\big\}<0.
\end{split}
\end{equation}
Then there exists $s\in (0,T)$ such that $$\max_{i\in \{1,\ldots,n_{\mathcal{X}}\}}\{g_i(\bm{\phi}_{\bm{x}_0,0}^{\bm{d}_{\bm{x}_0}}(s),s)\}=0,$$ implying that $$\bm{\phi}_{\bm{x}_0,0}^{\bm{d}_{\bm{x}_0}}(s)\in \partial \{\bm{x}\in \mathbb{R}^n\mid g_i(\bm{x},s)\leq 0,i=1,\ldots,n_{\mathcal{X}}\}$$ since $\partial \mathcal{X}_s=\partial \{\bm{x}\in \mathbb{R}^n\mid g_i(\bm{x},s)\leq 0,i=1,\ldots,n_{\mathcal{X}}\}=\cup_{i=1}^{n_{\mathcal{X}}}\{\bm{x}\in \mathcal{X}_s\mid g_i(\bm{x},s)=0\}$. Therefore, there exist $\bm{y}_1$ satisfying $g_j(\bm{y}_1,s)>0$, $j\in \{1,\ldots,n_{\mathcal{X}}\}$, and $\bm{x}_1$ satisfying $u(\bm{x}_1,0)\leq 0$ such that $\bm{\phi}_{\bm{x}_1,0}^{\bm{d}_{\bm{x}_0}}(s)=\bm{y}_1$, contradicting that all possible trajectories starting from the set $\{\bm{x}\in \mathbb{R}^n \mid u(\bm{x},0)\leq 0\}$ stay inside the set $\mathcal{X}_t$ at $t\in [0,T]$. Thus, $\bm{x}_0\in \partial \{\bm{x}\in \mathbb{R}^n \mid u(\bm{x},0)\leq 0\}$ and therefore $$\{\bm{x}\in \mathbb{R}^n\mid u(\bm{x},0)=0\}=\partial \{\bm{x}\in \mathbb{R}^n \mid u(\bm{x},0)\leq 0\}.$$

In conclusion, $\{\bm{x} \in \mathbb{R}^n\mid u(\bm{x},0)=0\}$ is the boundary of $\mathcal{R}_0=\{\bm{x}\in \mathbb{R}^n \mid u(\bm{x},0)\leq 0\}$. 
\hfill $\Box$
\end{IEEEproof}

Theorem \ref{convergence} states that the inner-approximation $\mathcal{X}_{0k}$ approximates the interior of the backward reachable set $\mathcal{R}_0$ with $k$ approaching infinity.
\begin{theorem}
\label{convergence}
Let the sequence of polynomials $\{\psi_{k+N}(\bm{x},t)\}_{k=1}^{\infty}$ satisfy Corollary \ref{unifrom}. Then the set $\mathcal{X}_{0k}$ satisfies that $\mathcal{X}_{0k}\subset \mathcal{R}_0$ and
 \[ \lim_{k\rightarrow \infty} \mu(\mathcal{R}_0\setminus \mathcal{X}_{0k})= \mu(\partial \mathcal{R}_0),\]
 where $\partial \mathcal{R}_0=\{\bm{x}\in \mathbb{R}^n\mid u(\bm{x},0)=0\}$, $\mathcal{X}_{0k}=\{\bm{x}\in \mathcal{X}_0\mid \psi_{k+N}(\bm{x},0)\leq 0\}$ and $u(\bm{x},t)$ is the viscosity solution to \eqref{HJB1}.
\end{theorem}
\begin{IEEEproof}
According to Corollary \ref{unifrom}, we have $\lim_{k\rightarrow +\infty}\int_{B(\bm{0},R)}|\psi_{k+N}(\bm{x},0)-u(\bm{x},0)|d\mu(\bm{x})=0$, implying that for every $\epsilon>0$,
$$\lim_{k\rightarrow +\infty}\mu(\{\bm{x}\mid |\psi_{N+k}(\bm{x},0)-u(\bm{x},0)|\geq \epsilon\})=0.$$ Then following the proof of Theorem 3 in \cite{lasserre2015} and combining with Lemma \ref{zero}, we have the conclusion. \hfill $\Box$
\end{IEEEproof}

\section{Examples and Discussions}
\label{examples}
In this section we evaluate our approach on three examples. All computations were performed on an i7-7500U 2.70GHz CPU with 32GB RAM running Windows 10. For numerical implementation, we formulate the sum-of-squares programming problem \eqref{sos} using the MATLAB package YALMIP \cite{Lofberg} and use Mosek \cite{mosek2015}\footnote{Mosek is free for academic use and can be obtained from {\tt https://www.mosek.com/}.} as a semi-definite programming solver. In order to evaluate the performance of our approach, we also present results for these three examples by dealing with \eqref{HJB1} directly. We employ the ROC-HJ solver \cite{bokanowski2013}\footnote{The ROC-HJ solver can be downloaded from \url{https://uma.ensta-paristech.fr/soft/ROC-HJ/}.} for solving \eqref{HJB1}.

\subsection{Examples}
\label{exs}
In this subsection we test our method on three illustrative examples. Examples \ref{ex1} and \ref{ex2} are employed to illustrate the performance of our method under different parameter settings. Example \ref{ex3} is primarily used to evaluate the scalability of our method. The parameters that control the performance of our approach in these three examples are presented in Table \ref{table}, which together shows the computation times for these three examples  in solving \eqref{HJB1} directly. Note that in solving \eqref{HJB1}, uniform grids of $500^2$ on the state space $[-1.1,1.1]\times [-1.1,1.1]$  are adopted for Examples \ref{ex1} and \ref{ex2}, and uniform grids of $10^7$ on the state space $[-0.55,0.55]^7$ for  Example \ref{ex3}. Due to the curse of dimensionality suffered by grid-based numerical methods for solving \eqref{HJB1}, this coarse grid for Example \ref{ex3} is adopted. 

\begin{table}[ht]
\begin{center}
\begin{tabular}{|c|c||c|c|c|c||c|}
  \hline
    \multirow{1}{*}{}&&\multicolumn{4}{|c||}{\texttt{SDP}~\eqref{sos}} & \multicolumn{1}{|c|}{\texttt{HJE}~\eqref{HJB1}} \\\hline
   Ex.&$[t_0,T]$&$k$&$d_{s}$& $d_{s'}$&  {Time}  & Time\\\hline
    \multirow{7}{*}{1} &\multirow{7}{*}{[0,1]} & &  & & \begin{tabular} {c|c}   1a  & 1b \end{tabular}&\multirow{7}{*}{334.20}\\
    \cline{3-6}
                                                                          & &4&2&2 &\begin{tabular} {c|c} 0.69&0.71 \end{tabular}&   \\
    \cline{3-6}
                                                                          & &6&4& 4&\begin{tabular} {c|c} 0.73& 0.83 \end{tabular}&   \\
     \cline{3-6}
                                                                          & &8&6& 6&\begin{tabular} {c|c} 0.81& 1.32 \end{tabular}&   \\
     \cline{3-6}
                                                                          & &10&8&8 &\begin{tabular} {c|c} 4.17 &4.08\end{tabular}&   \\
     \cline{3-6}
                                                                          & &12&10&10 &\begin{tabular} {c|c} 21.09& 21.20\end{tabular}&  \\
     \cline{3-6}
                                                                          & &14&12&12 &\begin{tabular} {c|c} 89.23 &97.43\end{tabular} & \\\hline
    \multirow{7}{*}{2}& \multirow{7}{*}{[0,1]} &   &   & & \begin{tabular} {c|c}   2a  & 2b\end{tabular}& \multirow{7}{*}{519.30}\\
       \cline{3-6}
                                                                         & &4&4& 4&\begin{tabular} {c|c} 0.60 &0.60\end{tabular}&     \\
       \cline{3-6}
                                                                         & &6&6& 6&\begin{tabular} {c|c} 1.07 &1.13\end{tabular}&   \\
       \cline{3-6}
                                                                         & &8&8& 8&\begin{tabular} {c|c} 3.80 &4.24\end{tabular}& \\
       \cline{3-6}
                                                                         & &10&10&10 &\begin{tabular} {c|c} 24.47 &25.90\end{tabular}&   \\
       \cline{3-6}
                                                                         &  &12&12& 12&\begin{tabular} {c|c} 88.54&93.11\end{tabular}& \\
       \cline{3-6}
                                                                         &  &14&14& 14&\begin{tabular} {c|c} 423.73&478.47\end{tabular}&\\ \hline
     \multirow{3}{*}{3}  & \multirow{3}{*}{[0,1]}&4&4&4 &~~~~57.45 &     \multirow{3}{*}{--} \\
            \cline{3-6}                                              & &5&4& 4&~~~~72.86 & \\
            \cline{3-6}
                                                                            & &6&6& 4&~~~~$5638.23$ & \\\hline
   \end{tabular}
\end{center}
\caption{\textit{Parameters and performance of our implementations of solving \eqref{sos} and \eqref{HJB1} on the examples presented in this section. $[t_0,T]$: reachability time interval; $k$: degree of the polynomial $\psi_k$ in \eqref{sos}; $d_s$: degree of sum-of-squares multipliers in the first constraint in \eqref{sos}; $d_{s'}$: degree of sum-of-squares multipliers in the second and third constraints in \eqref{sos}; Time: computation times in seconds.} }
\label{table}
\end{table}

\begin{example}
\label{ex1}
Consider a two-dimensional system given by
\begin{equation*}
\begin{aligned}
&\dot{x}=-\frac{1}{2}x-(\frac{1}{2}+d)y+\frac{1}{2},\\
&\dot{y}=-\frac{1}{2}y+1,
\end{aligned}
\end{equation*}
where $T=1$, $\mathtt{TR}=\{\bm{x}\mid  x^2+y^2-0.64 \leq 0\}$, $\mathcal{D}=\{d\mid 0.01^2-d^2\geq 0\}$, $\mathcal{X}_t=\{\bm{x}\mid x^2+y^2 -1 \leq 0\}$ for $t\in [0,T]$ and $1a). B(\bm{0},R)=\{\bm{x}\mid 1.21-(x^2+y^2)\geq 0\}$; $1b). B(\bm{0},R)=\{\bm{x}\mid 2-(x^2+y^2)\geq 0\}$.

\begin{figure*}[ht]
\centering
\setlength\fboxsep{0pt}
\setlength\fboxrule{0.15pt}
\begin{tabular}{cccc}
\fbox{\includegraphics[width=2.0in,height=2.0in]{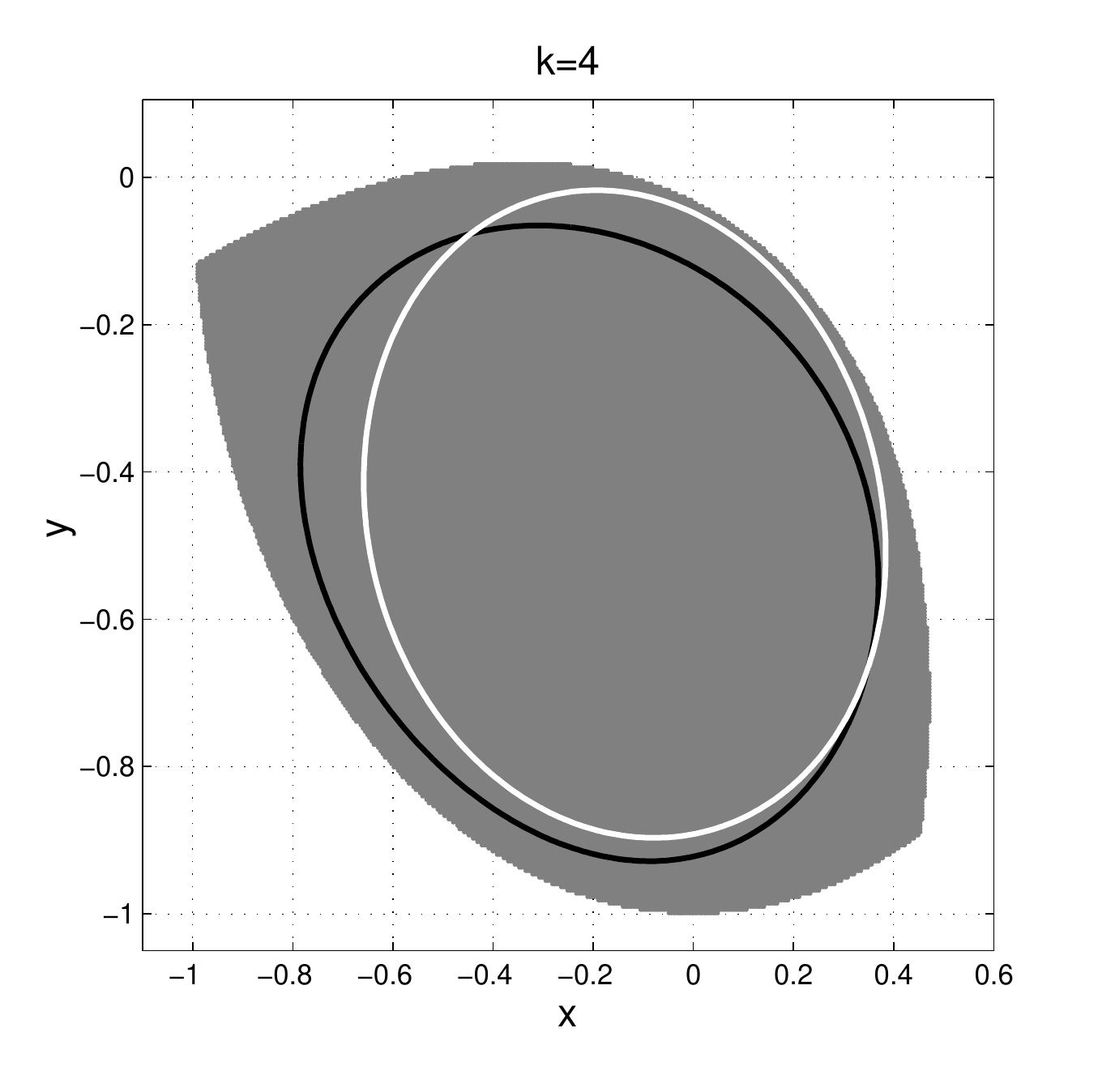}}
\fbox{\includegraphics[width=2.0in,height=2.0in]{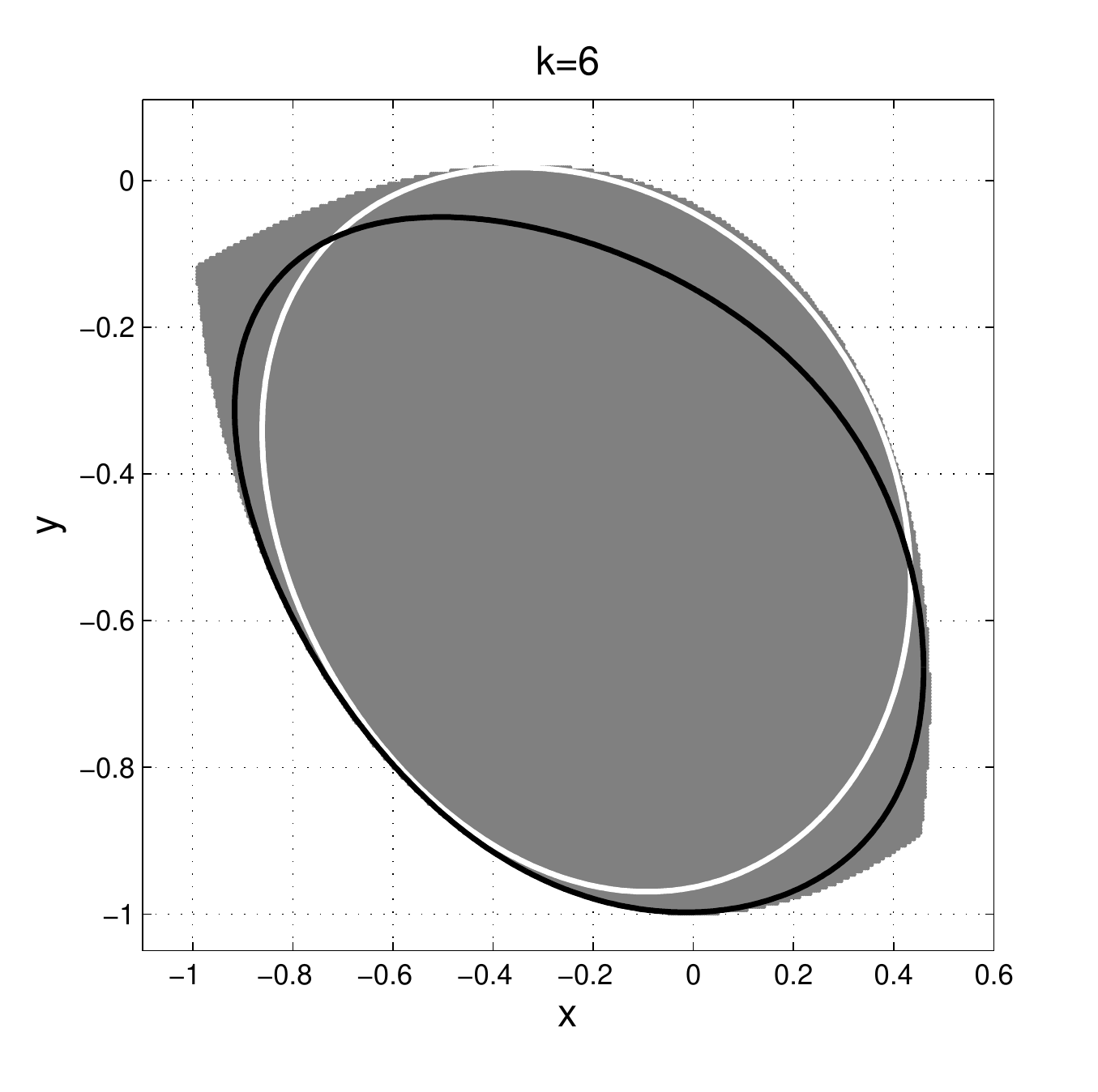}}
\fbox{\includegraphics[width=2.0in,height=2.0in]{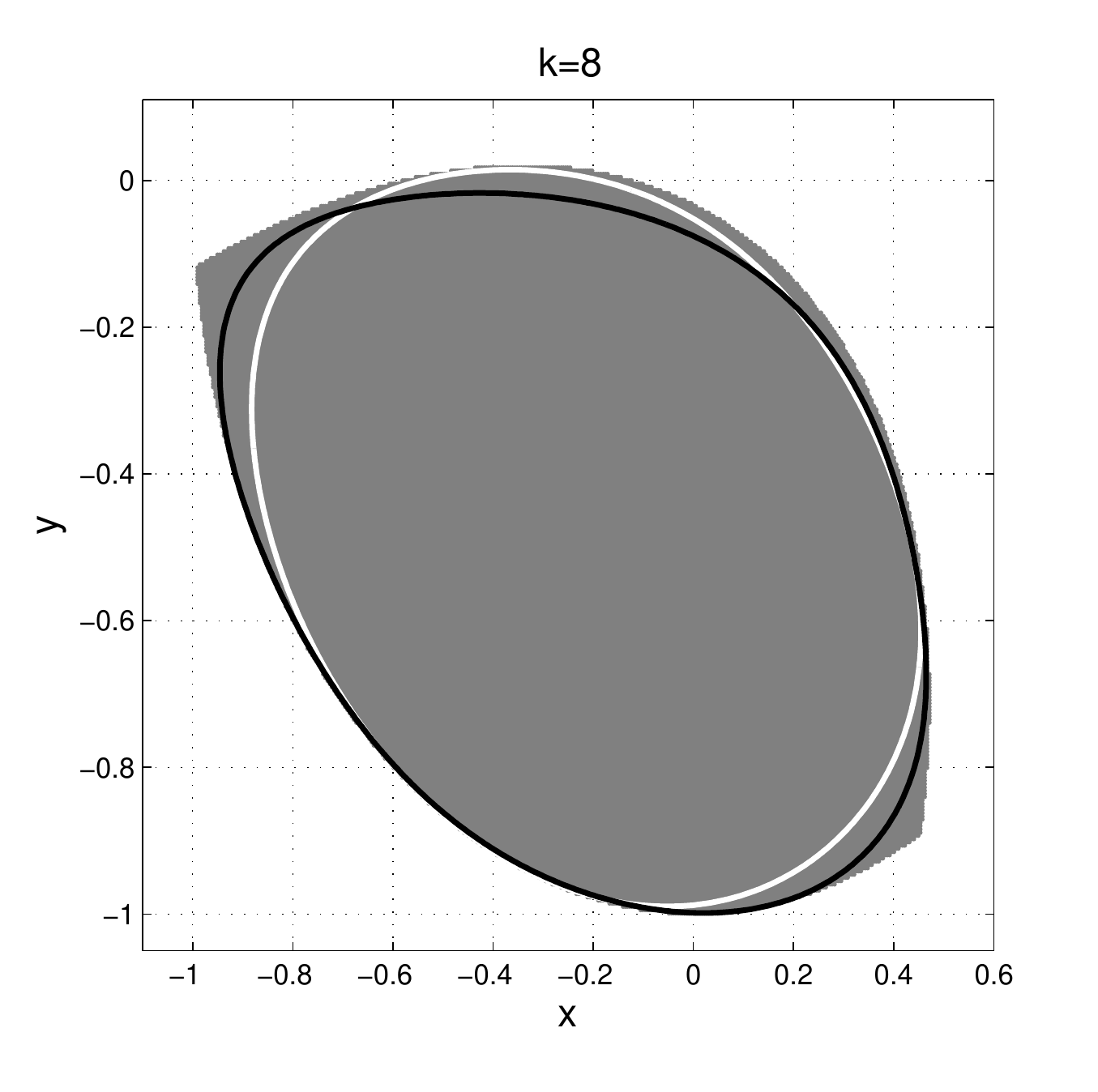}}\\
\fbox{\includegraphics[width=2.0in,height=2.0in]{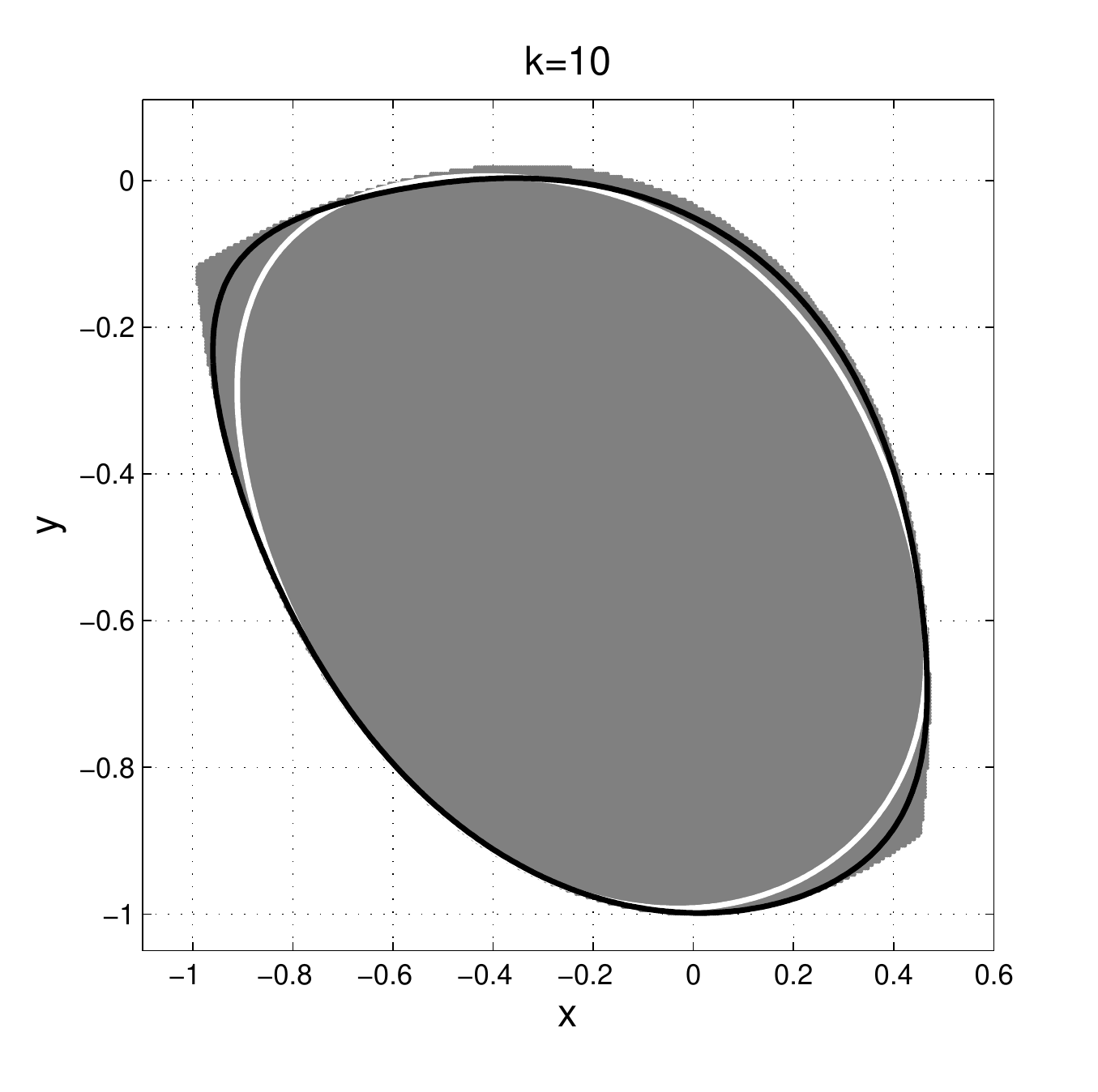}}
\fbox{\includegraphics[width=2.0in,height=2.0in]{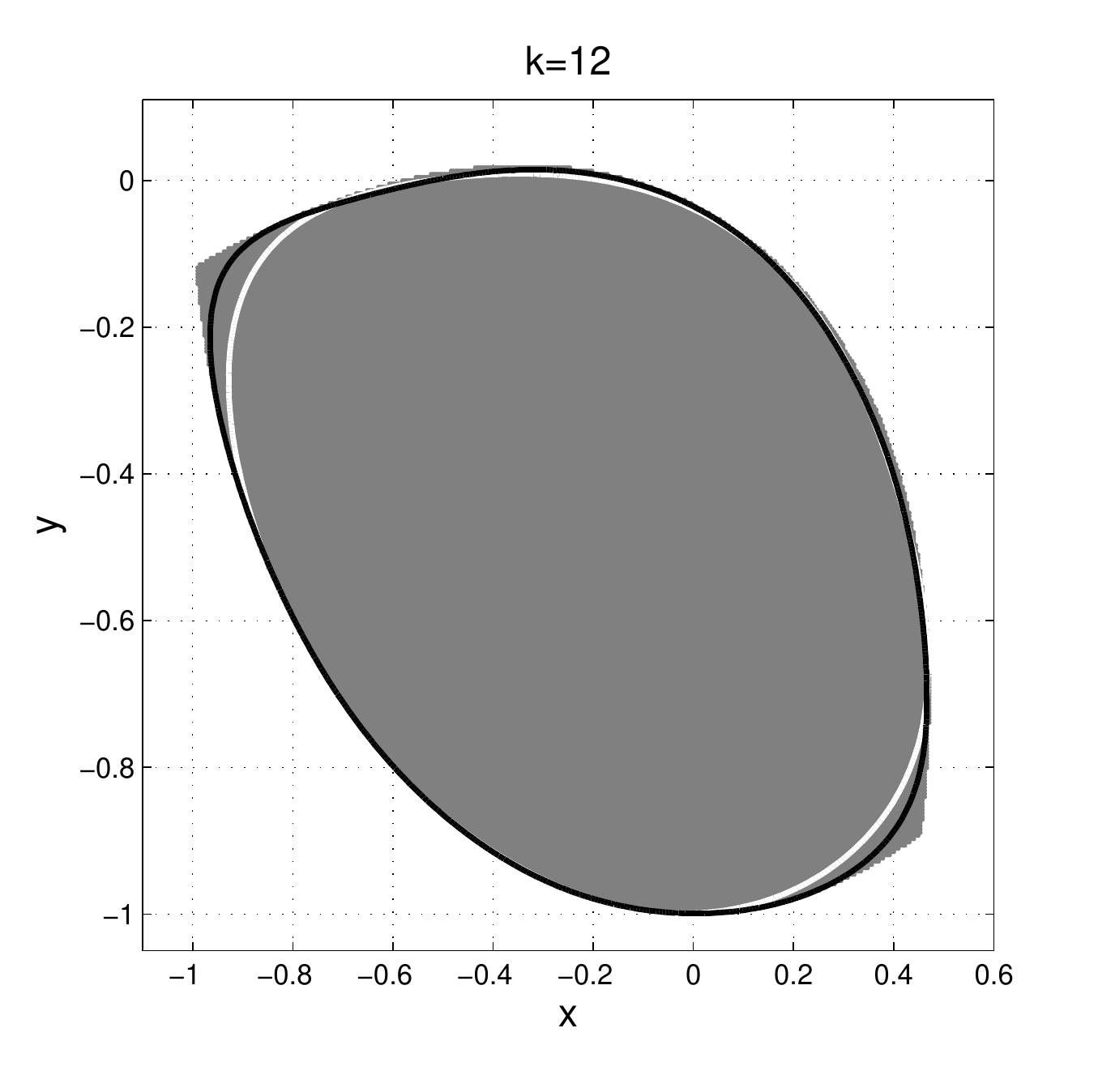}}
\fbox{\includegraphics[width=2.0in,height=2.0in]{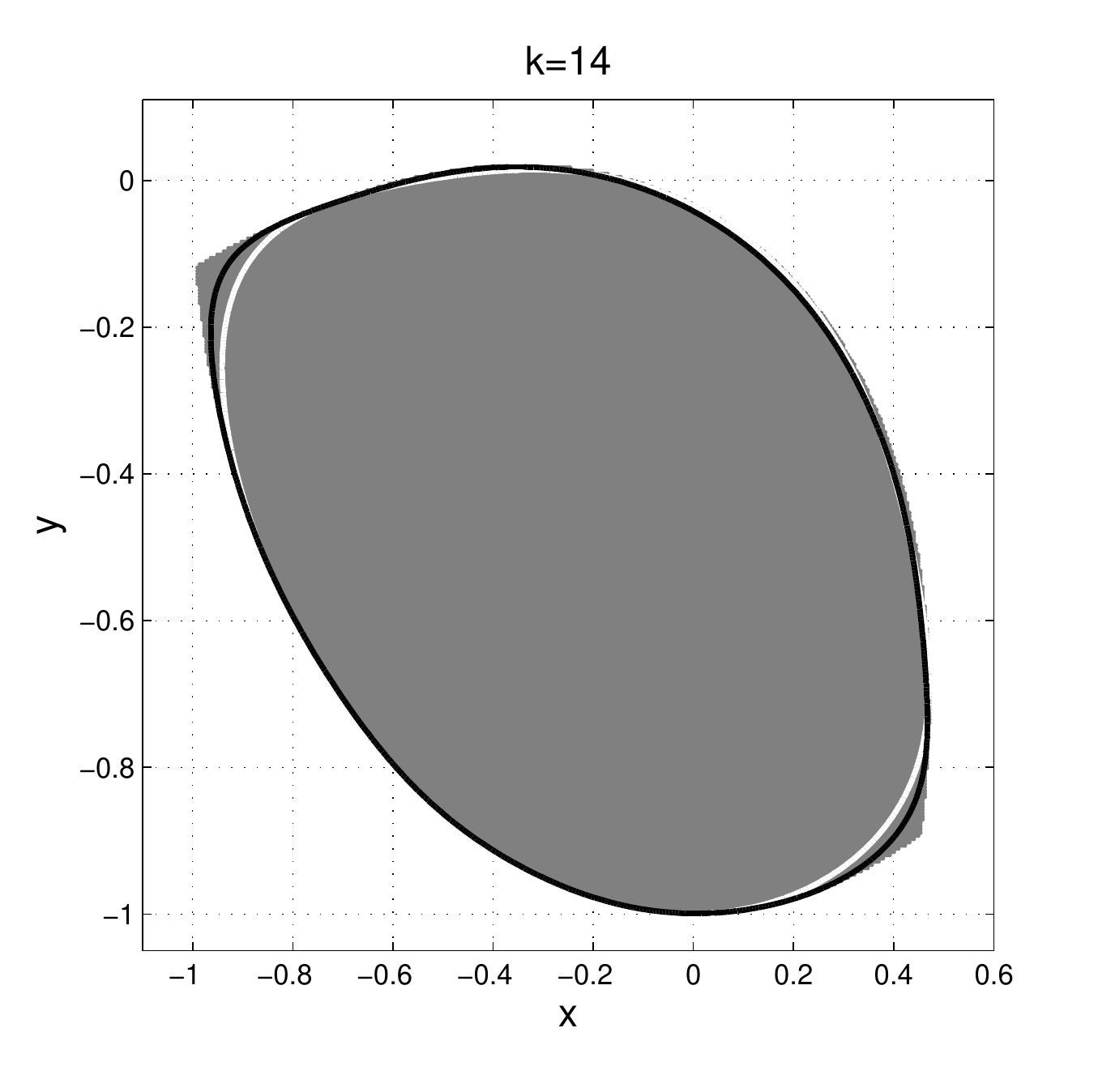}}
\end{tabular}
\caption{An illustration of computed backward reachable sets for Example \ref{ex1} at time $t=0$. \small{(Black and white curves denote the boundaries of inner-approximations in cases of 1a) and (1b), respectively. Gray region denotes the backward reachable set obtained via solving \eqref{HJB1}.)}}
\label{fig-one-1}
\end{figure*}
\end{example}

\begin{example}
\label{ex2}
Consider a scaled version of the reversed-time Van der Pol oscillator subject to uncertainties given by
\begin{equation*}
\begin{aligned}
&\dot{x}=-y,\\
&\dot{y}=0.4x+5(x^2-(d+0.2))y,
\end{aligned}
\end{equation*}
where $T=1$, $\mathtt{TR}=\{\bm{x}\mid  x^2+y^2- 0.25\leq 0\}$, $\mathcal{D}=\{d\mid 0.01^2-d^2\geq 0\}$, $\mathcal{X}_t=\{\bm{x}\mid x^2+y^2 - 0.64\leq 0\}$ for $t\in [0,T]$ and $2a). B(\bm{0},R)=\{\bm{x}\mid 0.65-(x^2+y^2)\geq 0\}$; $2b). B(\bm{0},R)=\{\bm{x}\mid 1-(x^2+y^2)\geq 0\}$.

\begin{figure*}[ht]
\centering
\setlength\fboxsep{0pt}
\setlength\fboxrule{0.15pt}
\begin{tabular}{ccc}
\fbox{\includegraphics[width=2.0in,height=2.0in]{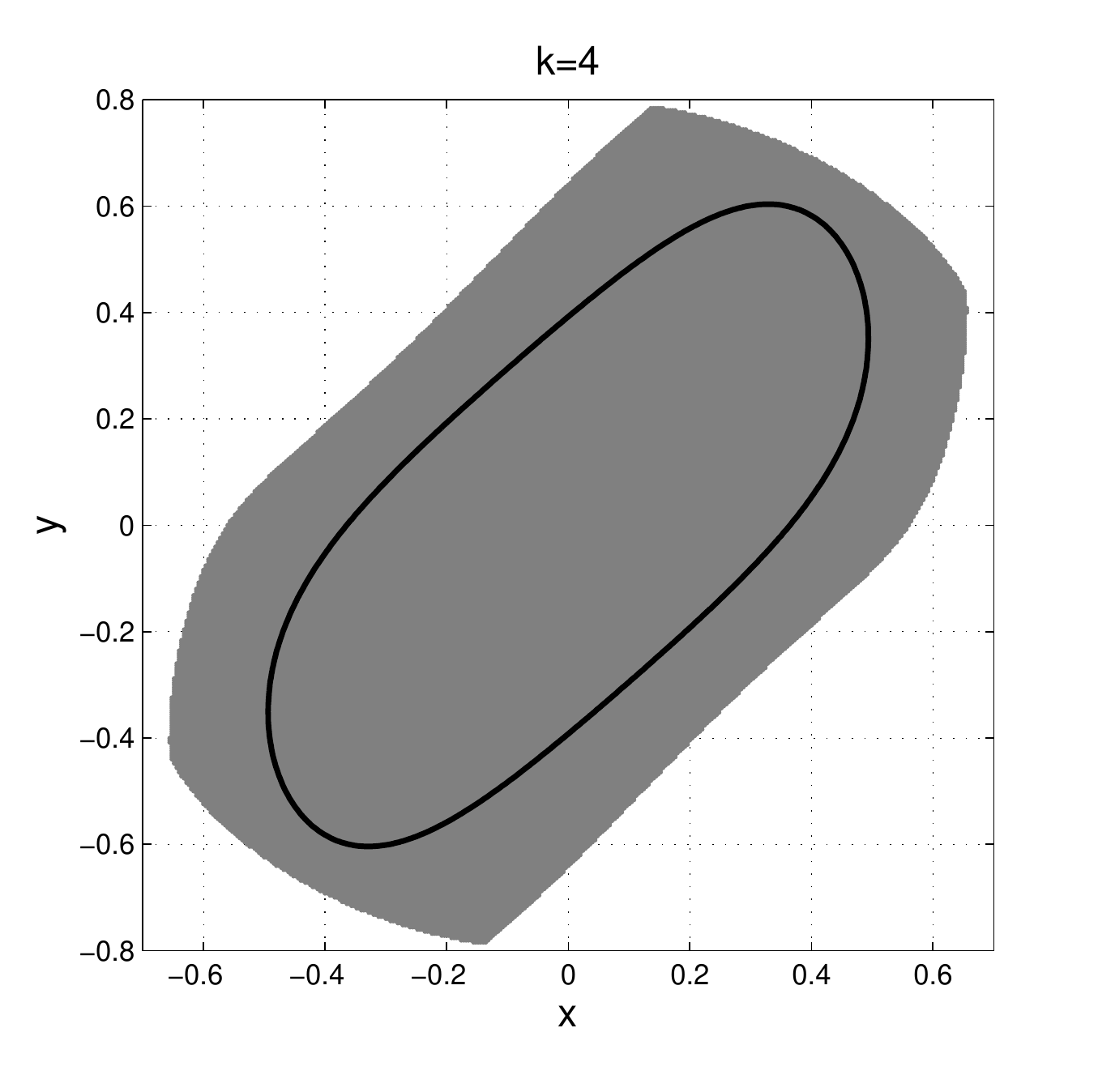}}
\fbox{\includegraphics[width=2.0in,height=2.0in]{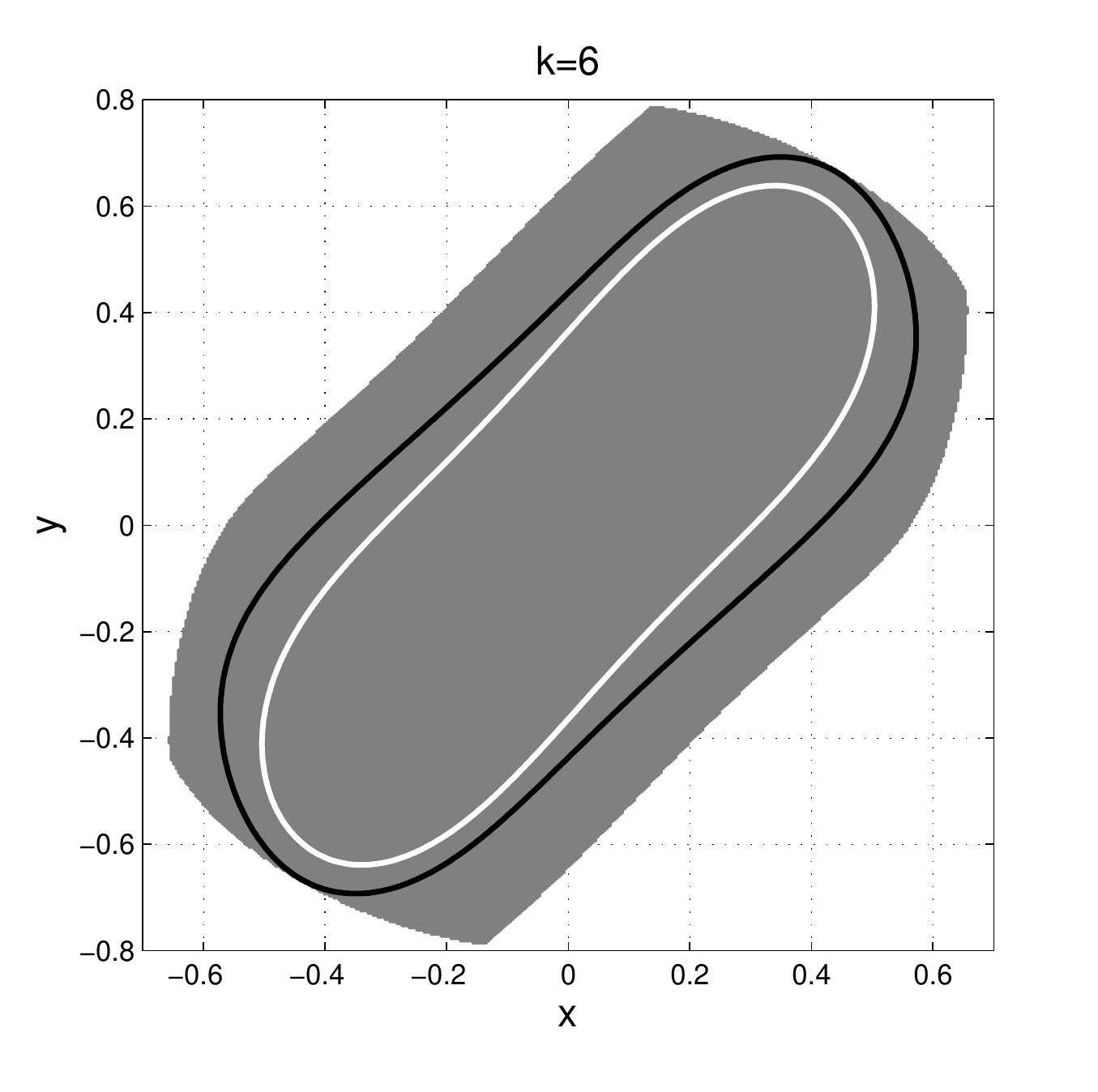}}
\fbox{\includegraphics[width=2.0in,height=2.0in]{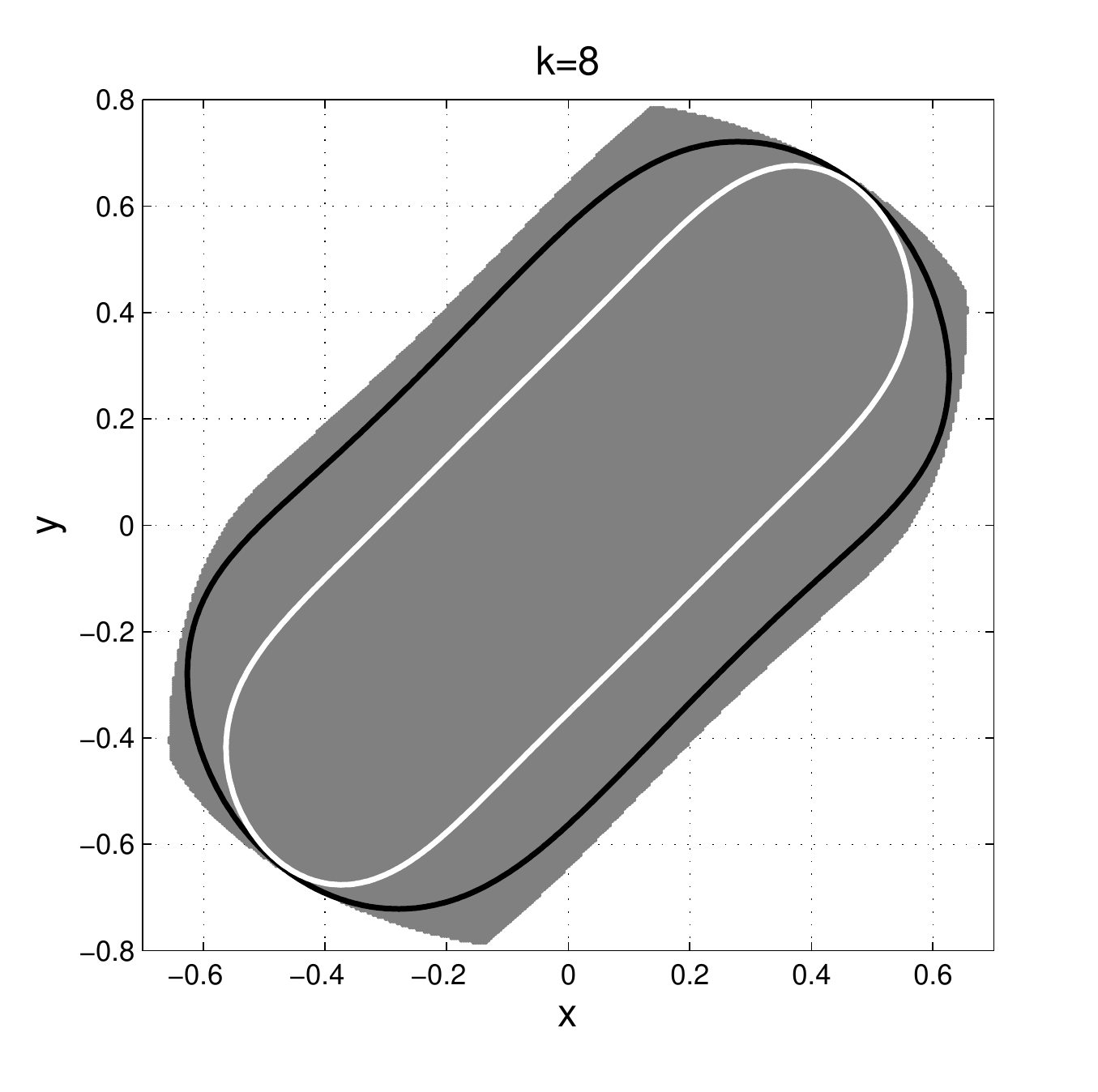}}\\
\fbox{\includegraphics[width=2.0in,height=2.0in]{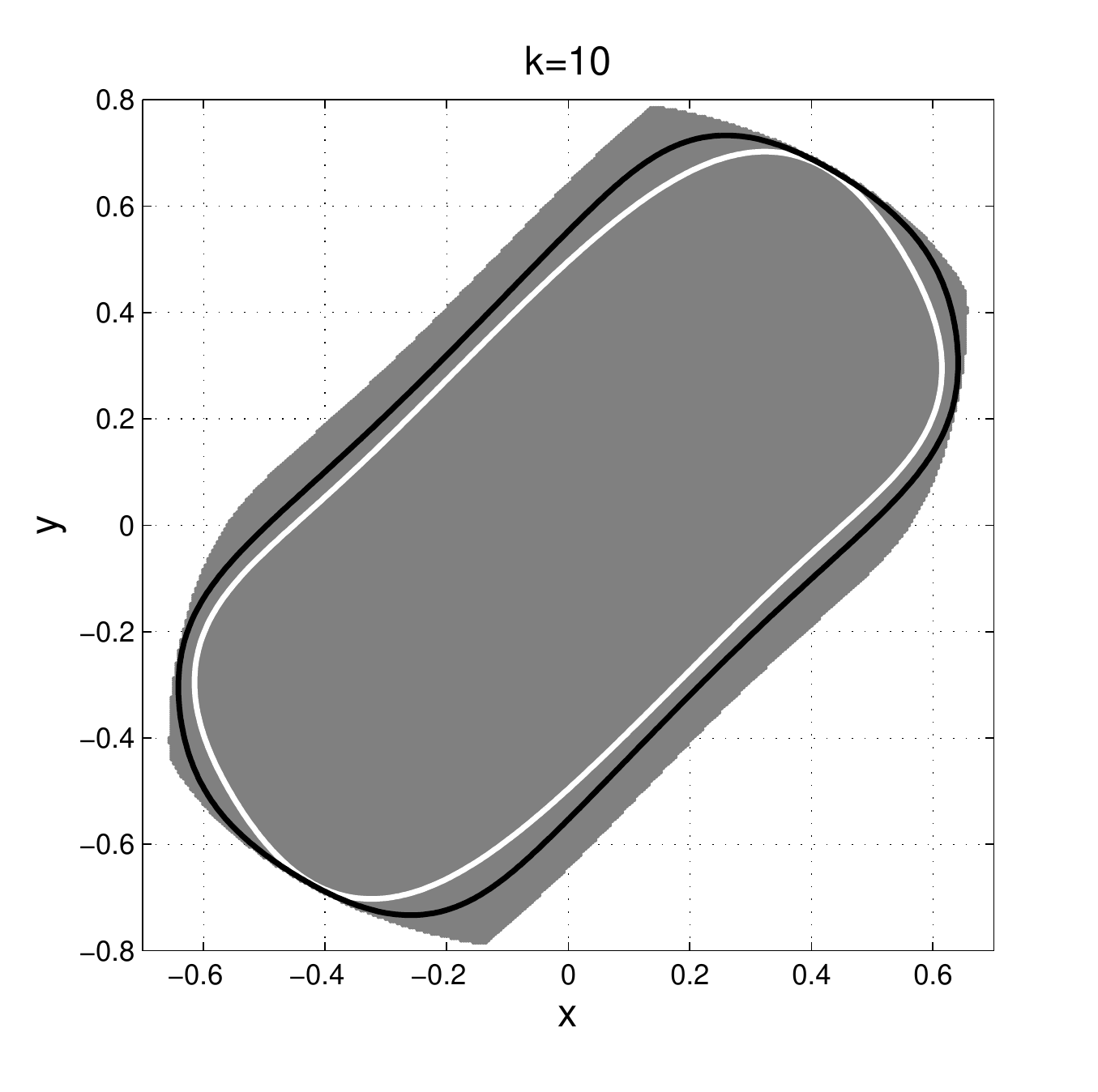}}
\fbox{\includegraphics[width=2.0in,height=2.0in]{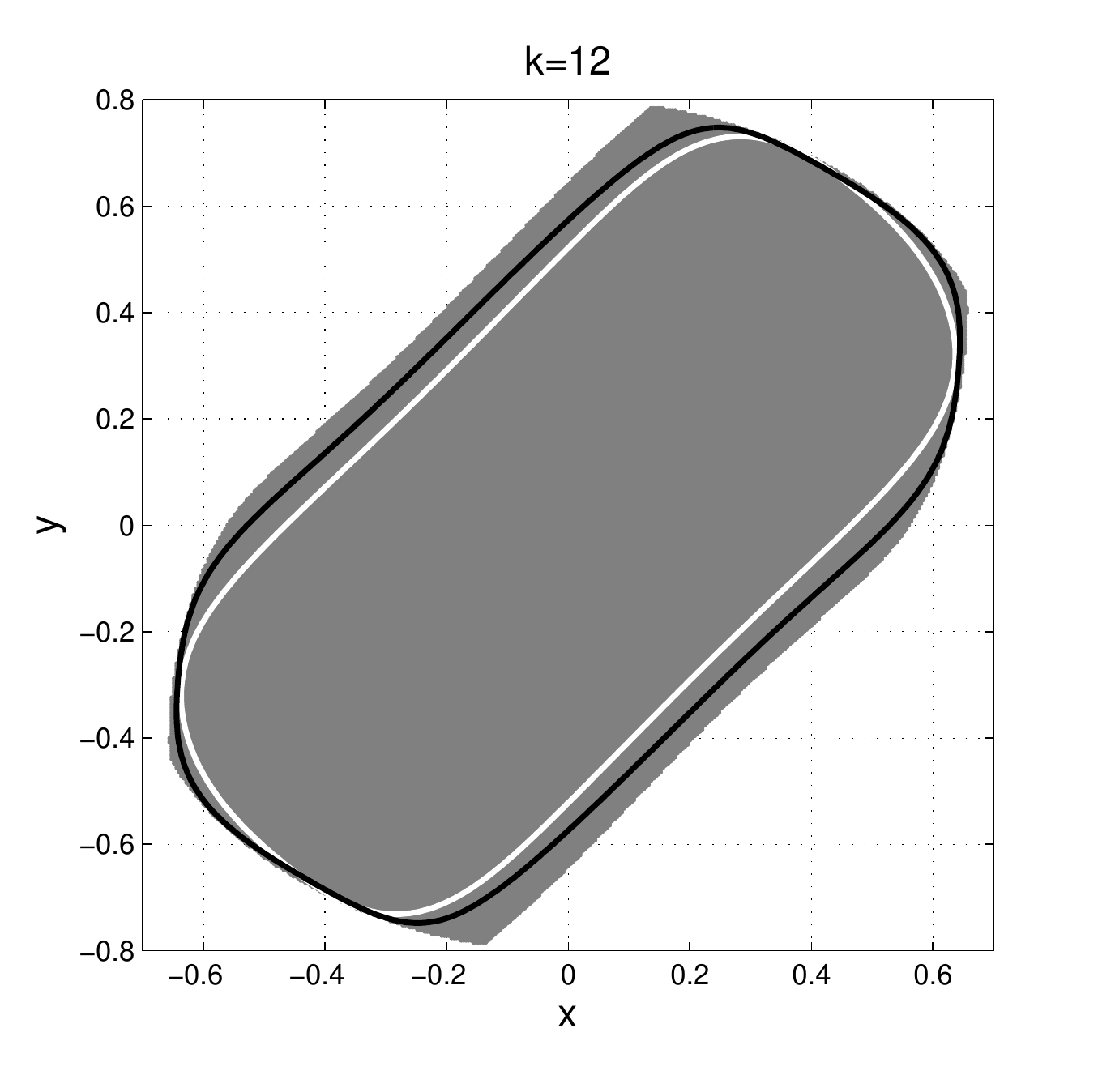}}
\fbox{\includegraphics[width=2.0in,height=2.0in]{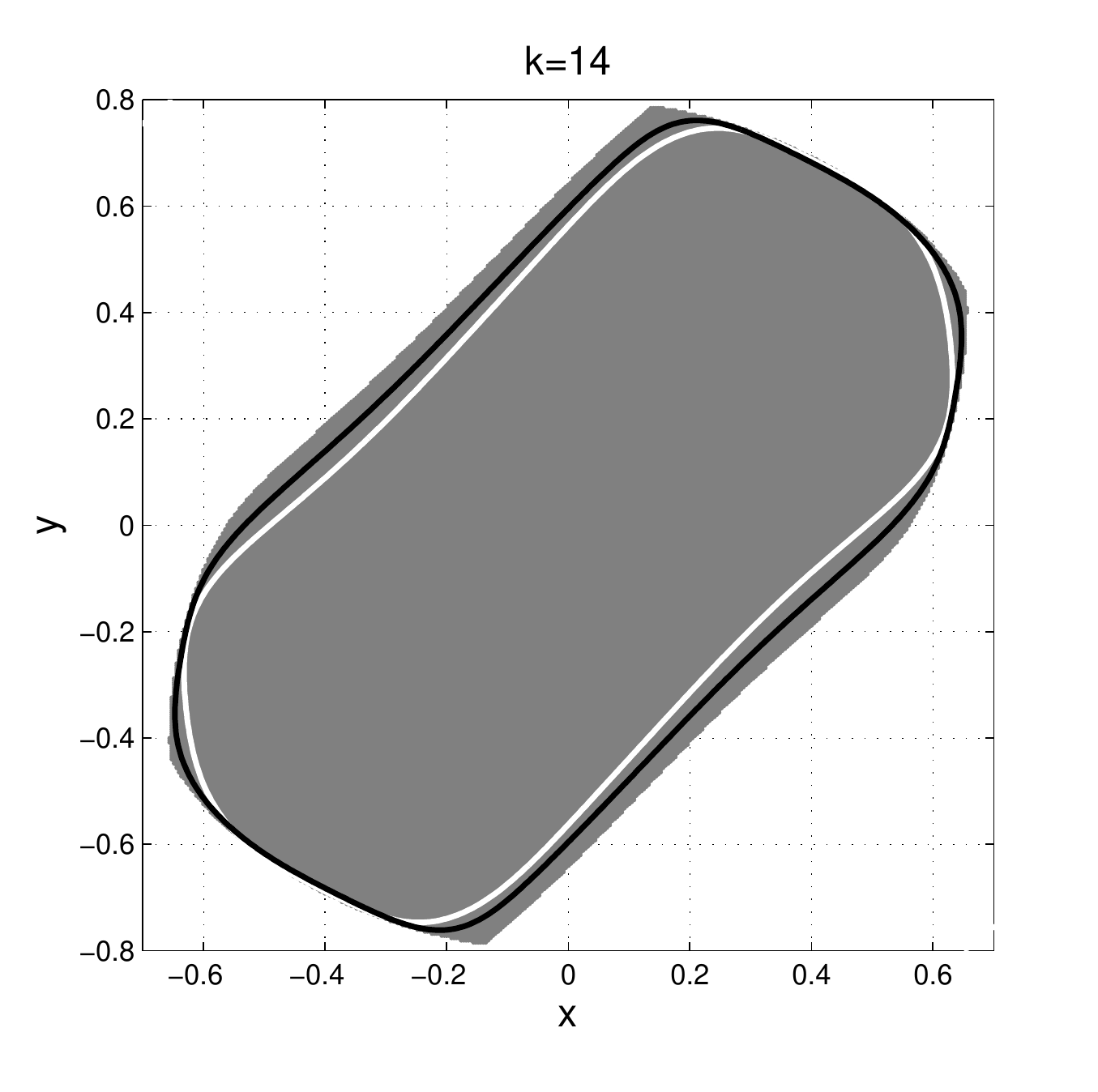}}
\end{tabular}
\caption{An illustration of computed backward reachable sets for Example \ref{ex1} at time $t=0$. \small{(Black and white curves denote the boundaries of inner-approximations in cases of 2a) and 2b), respectively. Gray region denotes the backward reachable set obtained via solving \eqref{HJB1}.)}}
\label{fig-one-2}
\end{figure*}
\end{example}

The computed inner approximations are illustrated in Fig. \ref{fig-one-1} and \ref{fig-one-2} for Examples \ref{ex1} and \ref{ex2} respectively. Note that the semi-definite program \eqref{sos} does not produce an inner approximation for Example \ref{ex2} when $k=4$ in the case of 2b) and consequently one cannot find the corresponding presentation in Fig. \ref{fig-one-2}. Observing the results illustrated in these two figures, we find that the accuracy of inner approximations to the backward reachable set is increasing with degree of the polynomial $\psi(\bm{x},t)$. Also, a relatively-fast convergence of inner-approximations to the backward reachable set is observed. The convergence rate is particularly fast when the degree of approximating polynomials is less than or equal to 10 and 12 for Example \ref{ex1} and Example \ref{ex2}, respectively. Moreover, the results in Fig. \ref{fig-one-2} indicate that  tighter sets $B(\bm{0},R)$ in \eqref{B} help to compute tighter inner approximations, although this indication is not obvious for Example \ref{ex1}.

Meanwhile, it is observed from Table \ref{table} that the semi-definite programming based method \eqref{sos} with polynomials of appropriate degree is more efficient in terms of computation time for Examples \ref{ex1} and \ref{ex2}, compared with grid-based numerical methods for solving \eqref{HJB1}. We continue exploring the performance of the semi-definite programming based method \eqref{sos} based on a seven-dimensional system.

\begin{example}
\label{ex3}
Consider an example adapted from a seven-dimensional biological system,
\begin{equation*}
\begin{aligned}
&\dot{x}_1=-0.4x_1+5x_3x_4+d,\\
&\dot{x}_2=0.4x_1-x2,\\
&\dot{x}_3=x_2-5x_3x_4,\\
&\dot{x}_4=5x_5x_6-5x_3x_4,\\
&\dot{x}_5=-5x_5x_6+5x_3x_4,\\
&\dot{x}_6=0.5x_7-5x_5x_6,\\
&\dot{x}_7=-0.5x_7+5x_5x_6,
\end{aligned}
\end{equation*}
where $T=1$, $\mathtt{TR}=\{\bm{x}\mid  x_1^2+(x_2+0.2)^2+x_3^2+x_4^2+x_5^2+x_6^2+x_7^2- 0.25 \leq 0\}$, $\mathcal{D}=\{d\mid 0.1^2-d^2\geq 0\}$, $\mathcal{X}_t=\{\bm{x}\mid  \sum_{i=1}^7x_i^2- 0.25 \leq 0\}$ for $t\in [0,T]$ and $B(\bm{0},R)=\{\bm{x}\mid  0.26- \sum_{i=1}^7x_i^2\geq 0\}$.

\begin{figure*}[ht]
\centering
\setlength\fboxsep{0pt}
\setlength\fboxrule{0.15pt}
\begin{tabular}{ccc}
\fbox{\includegraphics[width=3.5in,height=3.0in]{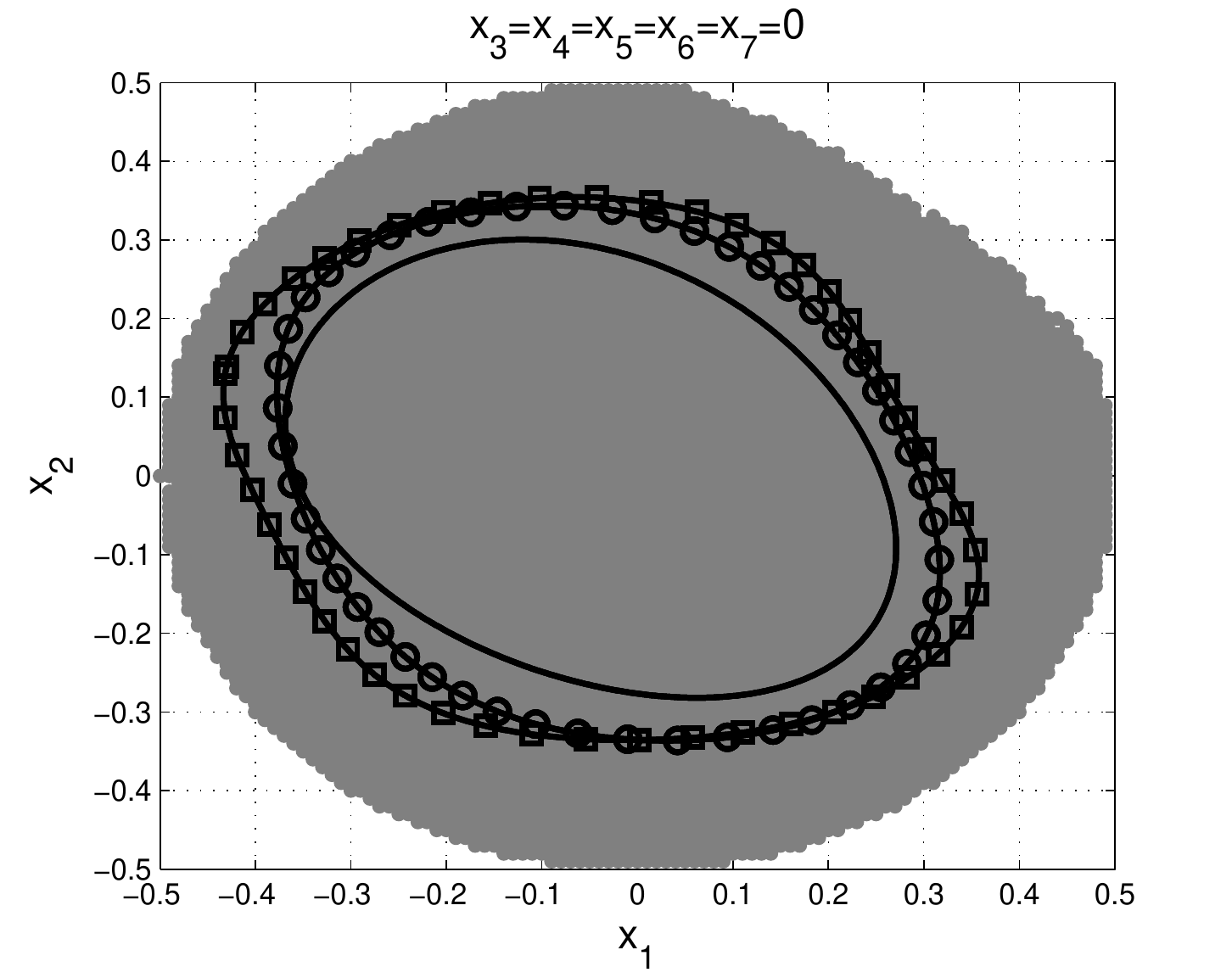}}
\fbox{\includegraphics[width=3.5in,height=3.0in]{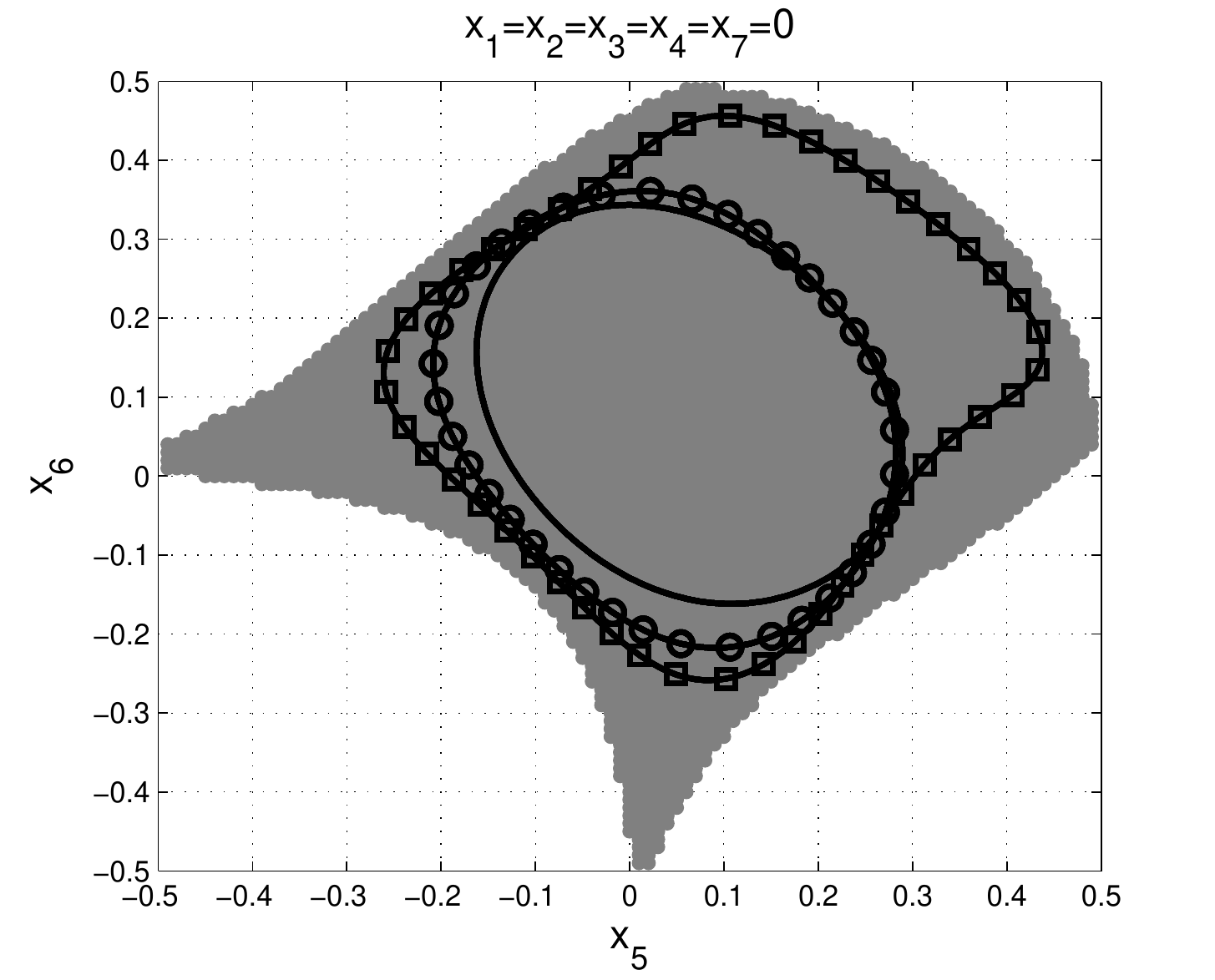}}
\end{tabular}
\caption{An illustration of computed backward reachable sets for Example \ref{ex3} at time $t=0$. \small{Black curve with square markers, black curve with circle markers and black curve denote the boundaries of inner-approximations when $k=6$, $k=5$ and $k=4$, respectively. Grey region denotes an estimation of the backward reachable set via simulation methods.)}}
\label{fig-two-3}
\end{figure*}
\end{example}

Unlike Examples \ref{ex1} and \ref{ex2}, the grid-based numerical method for solving \eqref{HJB1} runs out of memory and thus does not return an estimation for Example \ref{ex3}. In contrast, the method of solving \eqref{sos} is still able to compute inner-approximations, which are illustrated in Fig. \ref{fig-two-3}. Consequently, compared with grid-based numerical methods for solving \eqref{HJB1}, the semi-definite programming based method \eqref{sos} is capable of dealing with reachability problems of moderately high-dimensional systems, especially for cases where inner approximations formed by polynomials of low degree suffice. Although the size of the semidefinite program \eqref{sos} grows extremely
fast with the number of state and uncertainty variables and the degree of polynomials in \eqref{sos}, the computational efficiency and scalability advantage of the semi-definite programming based method \eqref{sos} could be further enhanced with template polynomials such as diagonally dominant sum-of-squares (DSOS) and scaled diagonally dominant-sum-of-squares (SDSOS) polynomials \cite{majumdar2014}.

Note that in order to shed light on the accuracy of computed inner approximations for Example \ref{ex3}, we partition state spaces $[-0.5,0.5]^2 \times [0,0]^5$ and $[0,0]^4 \times [-0.5,0.5]^2\times [0,0]$ and then employ the first-order Euler method to synthesize coarse estimations of the backward reachable set on planes $x_1-x_2$  with $x_3=x_4=x_5=x_6=x_7=0$ and $x_5-x_6$ with $x_1=x_2=x_3=x_4=x_7=0$ respectively. The estimations are the regions covered by grey points in Fig. \ref{fig-two-3}. The results in Fig. \ref{fig-two-3} further confirm that the accuracy of an inner approximation returned by solving \eqref{sos} is increasing with the degree of approximating polynomials. 

In Examples 1 to 3, we employ grid-based numerical methods for solving \eqref{HJB1}, e.g., Examples \ref{ex1} and \ref{ex2}, and simulation based methods, e.g., Example \ref{ex3}, to evaluate the quality of inner-approximations computed by solving \eqref{sos}. Another method is to estimate outer approximations of the backward reachable set $\mathcal{R}_0$ and calculate the Hausdorff distance between the outer approximation and inner approximation: narrower distance proves higher quality, which is the future work we are considering.

\section{Conclusion}
\label{conclusion}
We proposed a convex optimization based method to address the  problem of computing safe inner approximations of backward reachable sets for state-constrained polynomial systems subject to time-varying uncertainties in the setting of finite time horizons. The backward reachable set was first formulated as the zero sub-level set of the unique Lipschitz viscosity solution to a HJE. As opposed to traditionally grid-based numerical methods for solving the HJE, we proposed a novel semi-definite program, which was constructed from the HJE and falls within the convex programming category, to synthesize inner-approximations of the backward reachable set. We proved that solutions to the constructed semi-definite program are guaranteed to exist and can generate a convergent sequence of inner approximations to the interior of the backward reachable set. Three illustrative examples were employed to evaluate the performance of our approach.

\oomit{In this work we considered the reachability problem in the context of finite time horizons. However, reachability analysis over the infinite time horizon such as invariant generation (e.g., \cite{liu2011computing,topcu2010}) also has many important applications in practice such as safety controller design (e.g., \cite{FisacAZKGT17}). We therefore plan to extend the method in this paper to analyze the problem of infinite-time reachability in our future work.}

\textbf{Acknowledgements.}
We would like to thank the anonymous reviewers for their detailed and helpful comments, and Prof. Olivier Bokanowski for providing the ROC-HJ solver.
Bai Xue was funded partly by CAS Pioneer Hundred Talents
Program under grant No. Y8YC235015 and NSFC under
grant No. 61872341, Naijun Zhan was supported partly by
NSFC under grant No. 61625206 and 61732001, and Martin
Fr\"anzle was
supported by Deutsche Forschungsgemeinschaft under grant numbers DFG GRK
1765 (RTG SCARE) and FR 2715/4-1 (Science of Design of Societal-Scale
Cyber-Physical Systems).
\bibliographystyle{abbrv}
\bibliography{reference}
\section*{Appendix}
\textbf{The proof of Lemma \ref{p1}:}

\begin{proof}
Let $\bm{x}_1, \bm{x}_2 \in \mathcal{X}$, where $\mathcal{X}$ is an arbitrary but fixed compact set in $\mathbb{R}^n$, and $\bm{d}_1 \in \mathcal{M}_{t}$ such that
\begin{equation}
\begin{split}
u(\bm{x}_1,t)\leq &\max\big\{\max_{i\in \{1,\ldots,n_{\mathtt{TR}}\}}\{l_i(\bm{\phi}_{\bm{x}_1,t}^{\bm{d}_1}(T))\}, \max_{i\in \{1,\ldots,n_{\mathcal{X}}\}} \{\max_{s \in [t,T]}g_i(\bm{\phi}_{\bm{x}_1,t}^{\bm{d}_1}(s),s)\}\big\}+\epsilon,
\end{split}
\end{equation}
where $\epsilon>0$ is arbitrary but fixed.

Therefore, we infer that
\begin{equation}
\begin{split}
&u(\bm{x}_1,t)-u(\bm{x}_2,t)\\
&\leq \max\big\{\max_{i\in \{1,\ldots,n_{\mathtt{TR}}\}}\{l_i(\bm{\phi}_{\bm{x}_1,t}^{\bm{d}_1}(T))\}, \max_{i\in \{1,\ldots,n_{\mathcal{X}}\}} \{\max_{s \in [t,T]}g_i(\bm{\phi}_{\bm{x}_1,t}^{\bm{d}_1}(s),s)\}\big\}-u(\bm{x}_2,t)+\epsilon\\
&\leq \max\big\{\max_{i\in \{1,\ldots,n_{\mathtt{TR}}\}}\{l_i(\bm{\phi}_{\bm{x}_1,t}^{\bm{d}_1}(T))\}, \max_{i\in \{1,\ldots,n_{\mathcal{X}}\}} \{\max_{s \in [t,T]}g_i(\bm{\phi}_{\bm{x}_1,t}^{\bm{d}_1}(s),s)\}\big\}-\\
&~~~~\max\big\{\max_{i\in \{1,\ldots,n_{\mathtt{TR}}\}}\{l_i(\bm{\phi}_{\bm{x}_2,t}^{\bm{d}_1}(T))\},\max_{i\in \{1,\ldots,n_{\mathcal{X}}\}} \{\max_{s \in [t,T]}g_i(\bm{\phi}_{\bm{x}_2,t}^{\bm{d}_1}(s),s)\}\big\}+\epsilon \\
&\leq \max\big\{\max_{i\in \{1,\ldots,n_{\mathtt{TR}}\}}\{l_i(\bm{\phi}_{\bm{x}_1,t}^{\bm{d}_1}(T))-\l_i(\bm{\phi}_{\bm{x}_2,t}^{\bm{d}_1}(T))\}, \max_{i\in \{1,\ldots,n_{\mathcal{X}}\}}\{\max_{s \in [t,T]}(g_i(\bm{\phi}_{\bm{x}_1,t}^{\bm{d}_1}(s),s)-g_i(\bm{\phi}_{\bm{x}_2,t}^{\bm{d}_1}(s),s))\} \big\}+\epsilon\\
&\leq \max\big\{L_l\|\bm{\phi}_{\bm{x}_1,t}^{\bm{d}_1}(T)-\bm{\phi}_{\bm{x}_2,t}^{\bm{d}_1}(T)\|,  \max_{s\in [t,T]}L_g\|\bm{\phi}_{\bm{x}_1,t}^{\bm{d}_1}(s)-\bm{\phi}_{\bm{x}_2,t}^{\bm{d}_1}(s)\|\big\}+\epsilon\\
&\leq \max\{L_l,L_g\}e^{L_{\bm{F}}T}\|\bm{x}_1-\bm{x}_2\|+\epsilon,
\end{split}
\end{equation}
where $L_l$ and $L_g$ are Lipschitz constants \footnote{Such constants always exists since $l_i$ and $g_j$ are polynomials.} such that
\[|l_i(\bm{x}-\bm{y})|\leq L_l\|\bm{x}-\bm{y}\|, i=1,\ldots,n_{\mathtt{TR}}\]
\[|g_j(\bm{x}-\bm{y})|\leq L_g\|\bm{x}-\bm{y}\|, j=1,\ldots,n_{\mathcal{X}}\] over $\bm{x},\bm{y} \in S_{[0,T]}(\mathcal{X})$. $S_{[0,T]}(\mathcal{X})$ is defined in Definition \ref{re} and is compact. Note that we have used the following simple inequalities in above deduction:
\[\max\{A,B\}-\max\{C,D\}\leq \max\{A-C,B-D\}\] 
\[\max_s A(s)-\max_s B(s)\leq \max_s(A(s)-B(s)),\]
and the Gronwall-Bellman Lemma \cite{khalil2002}.

Using the similar argument as above with $\bm{x}_1$ and $\bm{x}_2$ reversed, we obtain
\[u(\bm{x}_2,t)-u(\bm{x}_1,t)\leq \max\{L_l,L_g\}e^{L_{F}T}\|\bm{x}_1-\bm{x}_2\|+\epsilon.\]
Since $\epsilon$ is arbitrary, we have
\begin{equation}
\label{x}
|u(\bm{x}_1,t)-u(\bm{x}_2,t)|\leq \max\{L_l,L_g\}e^{L_F T}\|\bm{x}_1-\bm{x}_2\|.
\end{equation}


Next, let $t_1,t_2\in [0,T]$, $t_1<t_2$ and $\bm{x}\in \mathcal{X}$, where $\mathcal{X}$ is an arbitrary but fixed compact set in $\mathbb{R}^n$. According to the definition of $u(\bm{x},t)$, i.e. \eqref{u},  there exists a $\bm{d}_1\in \mathcal{M}_t$ for any $\epsilon>0$ such that
\begin{equation}
\begin{split}
u(\bm{x},t_1)\leq &\max\big\{\max_{i\in \{1,\ldots,n_{\mathtt{TR}}\}}\{l_i(\bm{\phi}_{\bm{x},t_1}^{\bm{d}_1}(T))\}, \max_{i\in \in \{1,\ldots,n_{\mathcal{X}}\}}\{\max_{s\in [t_1,T]}g_{i}(\bm{\phi}_{\bm{x},t_1}^{\bm{d}_1}(s),s)\}\big\}+\epsilon.
\end{split}
\end{equation}
Therefore, we have the following deduction:
\begin{equation}
\begin{split}
&u(\bm{x},t_1)-u(\bm{x},t_2) \\
&\leq \max\{\max_{i\in \{1,\ldots,n_{\mathtt{TR}}\}}\{l_i(\bm{\phi}_{\bm{x},t_1}^{\bm{d}_1}(T))\}, \max_{i\in \{1,\ldots,n_{\mathcal{X}}\}}\{\max_{s\in [t_1,T]}g_i(\bm{\phi}_{\bm{x},t_1}^{\bm{d}_1}(s),s)\}\}-\\
&~~~~\max\{\max_{i\in \{1,\ldots,n_{\mathtt{TR}}\}}\{l_i(\bm{\phi}_{\bm{x},t_2}^{\bm{d}_1}(T))\},\max_{i\in \{1,\ldots,n_{\mathcal{X}}\}}\{\max_{s\in [t_2,T]}g_i(\bm{\phi}_{\bm{x},t_2}^{\bm{d}_1}(s),s)\}\}+\epsilon\\
&=\max\{\max_{i\in \{1,\ldots,n_{\mathtt{TR}}\}}\{l_i(\bm{\phi}_{\bm{x},t_1}^{\bm{d}_1}(T))\},\max_{i\in \{1,\ldots,n_{\mathcal{X}}\}}\{\max_{s\in [t_1,T]}g_i(\bm{\phi}_{\bm{x},t_1}^{\bm{d}_1}(s),s)\}\}-\\
&~~~~\max\{\max_{i\in \{1,\ldots,n_{\mathtt{TR}}\}}\{l_i(\bm{\phi}_{\bm{x},t_1}^{\bm{d}_1}(T-t_2+t_1))\},\max_{i\in \{1,\ldots,n_{\mathcal{X}}\}}\{\max_{s\in [t_1,T-t_2+t_1]}g_i(\bm{\phi}_{\bm{x},t_1}^{\bm{d}_1}(s),s)\}\}+\epsilon,
\end{split}
\end{equation}
where $L_l$ and $L_g$ are Lipschitz constants such that
\[|l_i(\bm{x})-l_i(\bm{y})|\leq L_l\|\bm{x}-\bm{y}\|, i=1,\ldots,n_{\mathtt{TR}}\]
\[|g_j(\bm{x})-g_j(\bm{y})|\leq L_g\|\bm{x}-\bm{y}\|, j=1,\ldots,n_{\mathcal{X}}\] over $\bm{x},\bm{y} \in S_{[0,T]}(\mathcal{X})$.

Assume that $s_0$ makes $\max_{i\in \{1,\ldots,n_{\mathcal{X}}\}}\{g_i(\bm{\phi}_{\bm{x},t_1}^{\bm{d}_1}(s_0),s_0)\}=\max_{i\in \{1,\ldots,n_{\mathcal{X}}\}}\{\max_{s\in [t_1,T]}g_i(\bm{\phi}_{\bm{x},t_1}^{\bm{d}_1}(s),s)\}$.

In case that $s_0\in [t_1,T-t_2+t_1]$, we have
\begin{equation}
\begin{split}
&u(\bm{x},t_1)-u(\bm{x},t_2)\\
&\leq \epsilon+\max\{\max_{i\in \{1,\ldots,n_{\mathtt{TR}}\}}\{l_i(\bm{\phi}_{\bm{x},t_1}^{\bm{d}_1}(T))\},\max_{i\in \{1,\ldots,n_{\mathcal{X}}\}}\{g_i(\bm{\phi}_{\bm{x},t_1}^{\bm{d}_1}(s_0),s_0)\}\}\\
&-\max\{\max_{i\in \{1,\ldots,n_{\mathtt{TR}}\}}\{l_i(\bm{\phi}_{\bm{x},t_1}^{\bm{d}_1}(T-t_1+t_2))\}, \max_{i\in \{1,\ldots,n_{\mathcal{X}}\}}\{g_i(\bm{\phi}_{\bm{x},t_1}^{\bm{d}_1}(s_0),s_0)\}\}\\
&\leq \max\big\{\max_{i\in \{1,\ldots,n_{\mathtt{TR}}\}}\{l_i(\bm{\phi}_{\bm{x},t_1}^{\bm{d}_1}(T))-l_i(\bm{\phi}_{\bm{x},t_1}^{\bm{d}_1}(T-t_2+t_1))\},0\big\}+\epsilon \\
&\leq \max\big\{L_l\|\bm{\phi}_{\bm{x},t_1}^{\bm{d}_1}(T)-\bm{\phi}_{\bm{x},t_1}^{\bm{d}_1}(T-t_2+t_1)\|,0\big\}+\epsilon\\
&\leq L_l b_{\bm{F}}|T-T+t_1-t_2|+\epsilon\\
&\leq L_l b_{\bm{F}}|t_2-t_1|+\epsilon,
\end{split}
\end{equation}
where $b_{\bm{F}}$ is an upper bound of $\|\bm{F}\|$ over $S_{[0,T]}(\mathcal{X})\times \mathcal{D}$.

In case that $s_0\in [T-t_2+t_1,T]$,  we first have
\begin{equation}
\begin{split}
&\max_{i\in \{1,\ldots,n_{\mathcal{X}}\}}\{\max_{s\in [t_1,T-t_2+t_1]}g_i(\bm{\phi}_{\bm{x},t_1}^{\bm{d}_1}(s),s)\}\geq \max_{i\in \{1,\ldots,n_{\mathcal{X}}\}}\{g_i(\bm{\phi}_{\bm{x},t_1}^{\bm{d}_1}(T-t_2+t_1),T-t_2+t_1)\}.
\end{split}
\end{equation}
Therefore, we have the following deduction:
\begin{equation}
\begin{split}
&u(\bm{x},t_1)-u(\bm{x},t_2)\\
&\leq \max\big\{\max_{i\in \{1,\ldots,n_{\mathtt{TR}}\}}\{l_i(\bm{\phi}_{\bm{x},t_1}^{\bm{d}_1}(T))\}, \max_{i\in \{1,\ldots,n_{\mathcal{X}}\}}\{g_i(\bm{\phi}_{\bm{x},t_1}^{\bm{d}_1}(s_0),s_0)\}\big\}-\\
&~~~~\max\big\{\max_{i\in \{1,\ldots,n_{\mathtt{TR}}\}}\{l_i(\bm{\phi}_{\bm{x},t_1}^{\bm{d}_1}(T-t_1+t_2))\}, \max_{i\in \{1,\ldots,n_{\mathcal{X}}\}}\{g_i(\bm{\phi}_{\bm{x},t_1}^{\bm{d}_1}(T-t_2+t_1),T-t_2+t_1)\}\big\}+\epsilon\\
&\leq \max\Big\{\max_{i\in \{1,\ldots,n_{\mathtt{TR}}\}}\big\{l_i(\bm{\phi}_{\bm{x},t_1}^{\bm{d}_1}(T))-l_i(\bm{\phi}_{\bm{x},t_1}^{\bm{d}_1}(T-t_2+t_1))\big\},\\
&~~~~\max_{i\in \{1,\ldots,n_{\mathcal{X}}\}}\big\{g_i(\bm{\phi}_{\bm{x},t_1}^{\bm{d}_1}(s_0),s_0)-g_i(\bm{\phi}_{\bm{x},t_1}^{\bm{d}_1}(T-t_2+t_1),T-t_2+t_1)\big\}\Big\}+\epsilon \\
&\leq \max\big\{L_l\|\bm{\phi}_{\bm{x},t_1}^{\bm{d}_1}(T)-\bm{\phi}_{\bm{x},t_1}^{\bm{d}_1}(T-t_2+t_1)\|, L_g\|\bm{\phi}_{\bm{x},t_1}^{\bm{d}_1}(s_0)-\bm{\phi}_{\bm{x},t_1}^{\bm{d}_1}(T-t_1+t_2)\| \big\}+\epsilon\\
&\leq \max\{L_l,L_g\} b_{\bm{F}}|T-T+t_2-t_1|+\epsilon\\
&\leq \max\{L_l,L_g\} b_{\bm{F}}|t_2-t_1|+\epsilon,
\end{split}
\end{equation}
where $b_{\bm{F}}$ is an upper bound of $\|\bm{F}\|$ over $S_{[0,T]}(\mathcal{X})\times \mathcal{D}$.

Therefore, $$u(\bm{x},t_1)-u(\bm{x},t_2)\leq \max\{L_l,L_g\} b_{\bm{F}}|t_2-t_1|+\epsilon$$ holds. Using the similar argument, we obtain that
\begin{equation}
u(\bm{x},t_1)-u(\bm{x},t_2)\geq -\max\{L_l,L_g\} b_{\bm{F}}|t_2-t_1|-\epsilon.
\end{equation}
Therefore,
\begin{equation}
\label{t}
|u(\bm{x},t_1)-u(\bm{x},t_2)|\leq \max\{L_l,L_g\} b_{\bm{F}}|t_2-t_1|
\end{equation}
since $\epsilon$ is arbitrary.



Combining inequalities \eqref{x} and \eqref{t}, we get:
\[|u(\bm{x}_1,t_1)-u(\bm{x}_2,t_2)|\leq \overline{C}(|t_1-t_2|+\|\bm{x}_1-\bm{x}_2\|)\]
for some constant $\overline{C}\geq 0$. Therefore $u(\bm{x},t)$ is locally Lipschitz continuous on $\mathbb{R}^n\times [0,T]$. 
\end{proof}

\textbf{The proof of Lemma \ref{dp}:}
\begin{proof}
Let
\begin{equation}
\label{v}
\begin{split}
w(\bm{x},t):=\sup_{\bm{d}\in \mathcal{M}_{[t,t+\delta]}}
&\max\big\{u(\bm{\phi}_{\bm{x},t}^{\bm{d}}(t+\delta),t+\delta),\max_{i\in \{1,\ldots,n_{\mathcal{X}}\}}\{\max_{s\in [t,t+\delta]}g_i(\bm{\phi}_{\bm{x},t}^{\bm{d}}(s),s)\}\big\}.
\end{split}
\end{equation}
We will prove that for $\forall \epsilon>0$, $|w-u|\leq \epsilon$.

According to the definition of $u(\bm{x},t)$, i.e. \eqref{u}, for any $\epsilon_1$, there exists $\bm{d}\in \mathcal{M}_t$ such that
\begin{equation}
\begin{split}
u(\bm{x},t)\leq \max\big\{&\max_{i\in \{1,\ldots,n_{\mathtt{TR}}\}}\{l_i(\bm{\phi}_{\bm{x},t}^{\bm{d}}(T))\},\max_{i\in \{1,\ldots,n_{\mathcal{X}}\}} \{\max_{s \in [t,T]}g_i(\bm{\phi}_{\bm{x},t}^{\bm{d}}(s),s)\}\big\}+\epsilon_1.
\end{split}
\end{equation}

We then separately define
$\bm{d}_1(s)$ as the restriction of $\bm{d}(s)\in \mathcal{M}_{t}$ over $[t,t+\delta]$ and $\bm{d}_2(s)$ as the restriction of $\bm{d}(s)\in \mathcal{M}_{t+\delta}$ over $s\in [t+\delta,T]$, and $\bm{y}=\bm{\phi}_{\bm{x},t}^{\bm{d}_1}(t+\delta)$,
we obtain that
\begin{equation}
\begin{split}
w(\bm{x},t)\geq &\max \Big\{u(\bm{y},t+\delta),\max_{i\in \{1,\ldots,n_{\mathcal{X}}\}}\{\max_{s\in [t,t+\delta]}g_i(\bm{\phi}_{\bm{x},t}^{\bm{d}_1}(s),s)\}\Big\}\\
\geq &\max\Big\{\max\big\{\max_{i\in \{1,\ldots,n_{\mathtt{TR}}\}}\{l_i(\bm{\phi}_{\bm{y},t+\delta}^{\bm{d}_2}(T))\}, \max_{i\in \{1,\ldots,n_{\mathcal{X}}\}}\{\max_{s\in [t+\delta,T]}g_i(\bm{\phi}_{\bm{y},t+\delta}^{\bm{d}_2}(s),s)\}\big\}, \max_{i\in \{1,\ldots,n_{\mathcal{X}}\}}\{\max_{s\in [t,t+\delta]}g_i(\bm{\phi}_{\bm{x},t}^{\bm{d}_1}(s),s)\big\}\Big\}\\
\geq & \max \Big\{\max_{i\in \{1,\ldots,n_{\mathtt{TR}}\}}\{l_i(\bm{\phi}_{\bm{x},t}^{\bm{d}}(T))\}, \max_{i\in \{1,\ldots,n_{\mathcal{X}}\}}\{\max_{s\in [t,T]}g_i(\bm{\phi}_{\bm{x},t}^{\bm{d}}(s),s)\}\Big\}\\
\geq & u(\bm{x},t)-\epsilon_1.
\end{split}
\end{equation}
Therefore,
\begin{equation}
\label{1}
u(\bm{x},t)\leq w(\bm{x},t)+\epsilon_1.
\end{equation}

By the definition of $w(\bm{x},t)$, i.e. \eqref{v}, for any $\epsilon_1>0$, there exists $\bm{d}_1(\cdot)\in \mathcal{M}_{[t,t+\delta]}$ such that
\begin{equation}
\begin{split}
w(\bm{x},t)\leq &\max\big\{u(\bm{\phi}_{\bm{x},t}^{\bm{d}_1}(t+\delta),t+\delta), \max_{i\in \{1,\ldots,n_{\mathcal{X}}\}}\{\max_{s\in [t,t+\delta]}g_i(\bm{\phi}_{\bm{x},t}^{\bm{d}_1}(s),s)\}\big\}+\epsilon_1.
\end{split}
\end{equation}
Also, by the definition of $u(\bm{x},t)$, i.e. \eqref{u}, for any $\epsilon_1$, there exists $\bm{d}_2(\cdot)\in \mathcal{M}_{t+\delta}$ such that
\begin{equation}
\begin{split}
u(\bm{y},t+\delta)\leq &\max\big\{\max_{i\in \{1,\ldots,n_{\mathtt{TR}}\}}\{l_i(\bm{\phi}_{\bm{y},t+\delta}^{\bm{d}_2}(T))\}, \max_{i\in \{1,\ldots,n_{\mathcal{X}}\}}\{\max_{s\in [t+\delta,T]}g_i(\bm{\phi}_{\bm{y},t+\delta}^{\bm{d}_2}(s),s)\}\big\}+\epsilon_1,
\end{split}
\end{equation}
where $\bm{y}=\bm{\phi}_{\bm{x},t}^{\bm{d}_1}(t+\delta)$.
We can define
\begin{equation}
\bm{d}(s)=
\begin{cases}
\bm{d}_1(s), s\in [t,t+\delta)\\
\bm{d}_2(s), s\in [t+\delta,T]
\end{cases},
\end{equation}
Therefore, we infer that
\begin{equation}
\label{2}
\begin{split}
w(\bm{x},t)&\leq \max\Big\{\max_{i\in \{1,\ldots,n_{\mathtt{TR}}\}}\{l_i(\bm{\phi}_{\bm{x},t}^{\bm{d}}(T))\}, \max_{i\in \{1,\ldots,n_{\mathcal{X}}\}}\{\max_{s\in [t,T]}g_i(\bm{\phi}_{\bm{x},t}^{\bm{d}}(s),s)\}\Big\}+2\epsilon_1\\
&\leq u(\bm{x},t)+2\epsilon_1.
\end{split}
\end{equation}

Combining \eqref{1} and \eqref{2} together, we finally have $|u-w|\leq \epsilon=2\epsilon_1$, implying that $u(\bm{x},t)=w(\bm{x},t)$ since $\epsilon_1$ is arbitrary.
\end{proof}

\textbf{The proof of Theorem \ref{solution}:}
\begin{proof}
Firstly, applying the definition of $u(\bm{x},s)$, i.e. \eqref{u}, to the terminal condition when $t=T$, $u(\bm{x},T)$ satisfies the boundary condition \eqref{HJB1}, i.e. $u(\bm{x},T)=\max\big\{\max_{i\in\{1,\ldots,n_{\mathtt{TR}}\}}\{l_i(\bm{x})\},\max_{i\in \{1,\ldots,n_{\mathcal{X}}\}}\{g_i(\bm{x},T)\}\big\}$.

The continuity property of the function $u(\bm{x},t)$ is assured by Lemma \ref{p1}. According to Definition \ref{viscosity}, a continuous function is a viscosity solution if it is both a sub-solution and a super-solution, we will respectively prove that $u$ is a viscosity sub- and super-solution to \eqref{HJB1}.

Firstly we prove that $u(\bm{x},t)$ is a sub-solution to \eqref{HJB1}. Let $\psi\in C^{\infty}(\mathbb{R}^n\times [0,T])$ such that $u-\psi$ attains a local maximum at $(\bm{x}_0,t_0)$, where $t_0\in [0,T]$; without loss of generality, assume that this maximum is 0, i.e. $u(\bm{x}_0,t_0)=\psi(\bm{x}_0,t_0)$. Therefore, there exists a positive value $\overline{\delta}$ such that $u(\bm{x},t)-\psi(\bm{x},t)\leq 0$ for $ (\bm{x},t)$ satisfying $\|\bm{x}-\bm{x}_0\|\leq \overline{\delta}$ and $0\leq t-t_0\leq \overline{\delta}$. Suppose \eqref{sub} in Definition \ref{viscosity} is false. Then there definitely exists a positive number $\epsilon_1$ such that
\begin{equation}
\label{con1}
\max_{{i\in \{1,\ldots,n_{\mathcal{X}}\}}}\{g_i(\bm{x}_0,t_0)\}\leq \psi(\bm{x}_0,t_0)-\epsilon_1
\end{equation}
holds. Therefore, there exists a sufficiently small $\delta'>0$ with $\delta'\leq \overline{\delta}$ such that
for $(\bm{x},\tau)$ satisfying $\|\bm{x}-\bm{x}_0\|\leq \delta'$ and $0\leq \tau-t_0\leq \delta'$,
$\max_{i\in \{1,\ldots,n_{\mathcal{X}}\}}\{g_i(\bm{x},\tau)\}\leq \psi(\bm{x}_0,t_0)-\frac{\epsilon_1}{2}.$
Also, there exists a positive number $\epsilon_2$ such that
\begin{equation}
\label{con2}
\partial_t \psi(\bm{x}_0,t_0)+H(\bm{x}_0,\triangledown_{\bm{x}} \psi) \leq -\epsilon_2
\end{equation}
holds. Since $\|\bm{\phi}_{\bm{x}_0,t_0}^{\bm{d}}(\tau)-\bm{x}_0\|=\|\int_{t=0}^{\tau}\bm{F}(\bm{x}(t),\bm{d}(t))dt\|\leq M (\tau-t_0)$ for $\bm{d}\in \mathcal{M}_{t_0}$, where $\tau\in [t_0,T]$ and  $M$ is a positive number such that $M\geq \|\bm{F}(\bm{x},\bm{d})\|$ over $S_{[t_0,T]}(\bm{x}_0)\times \mathcal{D}$, there exists small enough $\delta>0$ with $\delta\leq \delta'$ such that $\|\bm{\phi}_{\bm{x}_0,t_0}^{\bm{d}}(\tau)-\bm{x}_0\|\leq \delta'$ for $ \tau \in [0,\delta]$ and $\bm{d}\in \mathcal{M}_{[t_0,t_0+\delta]}$. Integrating \eqref{con2} from $t_0$ to $t_0+\delta$, we have 
\[\psi(\bm{\phi}_{\bm{x}_0,t_0}^{\bm{d}}(t_0+\delta),t_0+\delta)-\psi(\bm{x}_0,t_0)\leq -\frac{\epsilon_2}{2}\delta, \forall \bm{d}\in \mathcal{M}_{[t_0,t_0+\delta]}.\]
Further, since $u-\psi$ attains a local maximum of $0$ at $(\bm{x}_0,t_0)$, we infer that
\[u(\bm{\phi}_{\bm{x}_0,t_0}^{\bm{d}}(t_0+\delta),t_0+\delta)\leq u(\bm{x}_0,t_0)-\frac{\epsilon_2}{2}\delta\]
holds. Therefore, according to the dynamic principle \eqref{dy} in Lemma \ref{dp}, we finally have
\begin{equation}
\begin{split}
u(\bm{x}_0,t_0)=&\sup_{\bm{d}\in\mathcal{M}_{[t_0,t_0+\delta]}}
\max\{u(\bm{\phi}_{\bm{x}_0,t_0}^{\bm{d}}(t_0+\delta),t_0+\delta), \max_{i\in \{1,\ldots,n_{\mathcal{X}}\}}\{\max_{s\in [t_0,t_0+\delta]}g_i(\bm{\phi}_{\bm{x}_0,t_0}^{\bm{d}}(s),s)\}\}\\
&\leq u(\bm{x}_0,t_0)-\min\{\frac{\epsilon_1}{2},\frac{\epsilon_2}{2}\delta\}
\end{split}
\end{equation}
which is a contradiction, since $\epsilon_1$, $\epsilon_2$ and $\delta$ are positive. Therefore, $u$ is a sub-solution to \eqref{HJB1}.

Next, we prove that $u$ is a viscosity super-solution to \eqref{HJB1} as well. Let $\psi\in C^{\infty}(\mathbb{R}^n\times [0,T])$ such that $u-\psi$ attains a local minimum at $(\bm{x}_0,t_0)$, where $t_0\in [0,T]$. Similarly, we assume that this minimum is $0$, i.e. $u(\bm{x}_0,t_0)=\psi(\bm{x}_0,t_0)$. Therefore, there exists a positive value $\overline{\delta}$ such that $u(\bm{x},t)-\psi(\bm{x},t)\geq 0$ for $ (\bm{x},t)$ satisfying $\|\bm{x}-\bm{x}_0\|\leq \overline{\delta}$ and $0\leq t-t_0\leq \overline{\delta}$.

If \eqref{super} is false, then either
\begin{equation}
\label{con21}
\max_{i\in \{1,\ldots,n_{\mathcal{X}}\}}\{g_i(\bm{x}_0,t_0)\}\geq \psi(\bm{x}_0,t_0)+\epsilon_1,
\end{equation}
holds
or
\begin{equation}
\label{con22}
\partial_t \psi(\bm{x}_0,t_0)+H(\bm{x}_0,\nabla_{\bm{x}}u)\geq \epsilon_2,
\end{equation}
holds for some $\epsilon_1$, $\epsilon_2>0$.

If \eqref{con21} holds, then there is a small enough $\delta'>0$ with $\delta'\leq \overline{\delta}$ such that for $(\bm{x},t)$ satisfying $t_0\leq t\leq t_0+\delta'$ and
$\|\bm{x}-\bm{x}_0\|\leq \delta'$,
\[\max_{{i\in \{1,\ldots,n_{\mathcal{X}}\}}}\{g_i(\bm{x},t)\}\geq \psi(\bm{x}_0,t_0)+\frac{\epsilon_1}{2}=u(\bm{x}_0,t_0)+\frac{\epsilon_1}{2}.\]
Moreover, there exists $\delta>0$ with $\delta\leq \delta'$ such that 
$\|\bm{\phi}_{\bm{x}_0,t_0}^{\bm{d}}(\tau)-\bm{x}_0\|\leq \delta'$ for $ \tau \in [0,\delta]$ and $\bm{d}\in \mathcal{M}_{[t_0,t_0+\delta]}$.

Then the dynamic programming principle \eqref{dy} in Lemma \ref{dp} yields
\begin{equation}
\begin{split}
u(\bm{x}_0,t_0)=&\sup_{\bm{d}\in \mathcal{M}_{[t_0,t_0+\delta]}}\max\{u(\bm{\phi}_{\bm{x}_0,t_0}^{\bm{d}}(t_0+\delta),t_0+\delta), \max_{i\in \{1,\ldots,n_{\mathcal{X}}\}}\{\max_{s\in [t_0,t_0+\delta]} g_i(\bm{\phi}^{\bm{d}}_{\bm{x}_0,t_0}(s),s)\}\}\\
&\geq u(\bm{x}_0,t_0)+\frac{\epsilon_1}{2},
\end{split}
\end{equation}
which is a contradiction since $\epsilon_1>0$.

However, if \eqref{con22} holds, there is a small enough $\delta'>0$ with $\delta'\leq \overline{\delta}$ such that $\partial_t \psi(\bm{x},t)+H(\bm{x},\nabla_{\bm{x}}u)\geq \frac{\epsilon_2}{2}$ for $ (\bm{x},t)$ satisfying $\|\bm{x}-\bm{x}_0\|\leq \delta'$ and $0\leq t-t_0\leq \delta'$. Moreover, there exists a postive $\delta\leq \delta'$ such that $\|\bm{\phi}_{\bm{x}_0,t_0}^{\bm{d}}(\tau)-\bm{x}_0\|\leq \delta'$ for $\tau \in [0,\delta]$ and $\bm{d}\in \mathcal{M}_{[t_0,t_0+\delta]}$. Therefore, there exists $\bm{d}_1\in \mathcal{M}_{[t_0,t_0+\delta]}$ such that
\begin{equation}
\label{exist}
\delta\frac{\epsilon_2}{4}\leq \psi(\bm{\phi}_{\bm{x}_0,t_0}^{\bm{d}_1}(t_0+\delta),t_0+\delta)-\psi(\bm{x}_0,t_0).
\end{equation}
Note that \eqref{exist} can be assured by integrating \eqref{con22} with the fixed strategy $\bm{d}_1$ from $t_0$ to $t_0+\delta$ since $\psi\in C^{\infty}(\mathbb{R}^n\times [0,T])$.
Further, due to the fact that $u-\psi$ attains a local minimum at $(\bm{x}_0,t_0)$ and $u(\bm{x}_0,t_0)=\psi(\bm{x}_0,t_0)$, we have
\begin{equation}
\delta\frac{\epsilon_2}{4}\leq u(\bm{\phi}_{\bm{x}_0,t_0}^{\bm{d}_1}(t_0+\delta),t_0+\delta)-u(\bm{x}_0,t_0).
\end{equation}
Therefore, the following contradiction is obtained with the help of \eqref{dy}:
\begin{equation}
\begin{split}
u(\bm{x}_0,t_0)=&\sup_{\bm{d}\in \mathcal{M}_{[t_0,t_0+\delta]}}\{\max(u(\bm{\phi}_{\bm{x}_0,t_0}^{\bm{d}}(t_0+\delta),t_0+\delta),  \max_{i\in \{1,\ldots,n_{\mathcal{X}}\}}\{\max_{s\in[t_0,t_0+\delta]}g_i(\bm{\phi}_{\bm{x}_0,t_0}^{\bm{d}}(s),s))\}\}\\
&\geq u(\bm{x}_0,t_0) + \frac{\epsilon_2}{4}\delta,
\end{split}
\end{equation}
Thus, \eqref{super} holds and $u$ is a super-solution to \eqref{HJB1}. 

Uniqueness follows from Proposition 2 in \cite{fialho1999}: Let $v(\bm{x},t):\mathbb{R}^n\times [0,T]\mapsto \mathbb{R}$ and $u(\bm{x},t): \mathbb{R}^n\times [0,T]\mapsto \mathbb{R}$ be viscosity solutions with identical boundary condition $v(\bm{x},T)=u(\bm{x},T)$.
According to local comparison principle as illustrated in Proposition 1 in \cite{fialho1999}, $v\leq u$ and $u\leq v$, implying that $u=v$ over $\mathbb{R}^n\times [0,T]$.
\end{proof}

\end{document}